\documentclass[12pt]{article}
\usepackage{amsmath,amssymb}
\newtheorem{theorem}{Theorem}[section]
\newtheorem{lemma}{Lemma}[section]
\newtheorem{definition}{Definition}[section]

\newtheorem{prop}{Proposition}

\newcommand{\eps}{\epsilon}

\newcommand{\G}{\Gamma}

\newcommand{{\z}}{\mathbb Z}
\newcommand{\R}{\mathbb R}
\newcommand{\N}{\mathbb N}

\newcommand{\Pro}{\mathbb P}
\newcommand{\Z}{\mathbb Z}
\newcommand{\F}{\mathbb F}
\newcommand{\E}{\mathbb E}

\newcommand{\varcon}{\tilde{c}}

\newcommand{\fbm}{S}
\newcommand{\walk}{X}
\newcommand{\kmac}{K}

\long\def\note#1/{{\bf [#1]} }
\title{Power law P{\'o}lya's urn and fractional Brownian motion}
\author{Alan Hammond\thanks{University of Oxford.
Research undertaken while at Courant Institute of Mathematical Sciences and 
supported in part by the US National Science Foundation under
grants OISE-07-30136 and DMS-0806180.}
\ and Scott Sheffield\thanks{Courant Institute of Mathematical Sciences and
  Massachusetts Institute of Technology. Partially supported
by NSF grants DMS-0403182, DMS-0645585 and OISE-07-30136.}}

\date{}

\begin{document}
\maketitle
\begin{abstract} 
We introduce a natural family of random walks $S_n$ on $\mathbb Z$ that
scale to fractional Brownian motion.  The increments $\walk_n := S_n - S_{n-1}
\in \{\pm 1\}$ have the property that given $\{ \walk_k : k < n \}$, the
conditional law of $\walk_n$ is that of $\walk_{n - k_n}$, where $k_n$ is sampled
independently from a fixed law $\mu$ on the positive integers.  When $\mu$
has a roughly power law decay (precisely, when $\mu$ lies in the domain of
attraction of an $\alpha$-stable subordinator, for $0<\alpha<1/2$) the walks
scale to fractional Brownian motion with Hurst parameter $\alpha + 1/2$.
The walks are easy to simulate and their increments satisfy an FKG
inequality.  In a sense we describe, they are the natural ``fractional''
analogues of simple random walk on $\mathbb Z$.
\end{abstract}

%\begin{document}
%\begin{center}
%{Weighted spanning forests on the integers and  fractional Brownian motion}\bs\bs\\
%{Alan Hammond and Scott Sheffield}\bs\bs
%\end{center}

\begin{section}{Introduction}

Fractional Brownian motion is a one-parameter family of stochastic
processes, mapping the real line to itself, that are the only 
stationary-increment
Gaussian processes that, for some fixed $H
>0$, are invariant under the space-time rescalings $\fbm \to Y$ of the
form $Y_t = c^{-H} \fbm_{c t}$ with $c>0$. The parameter $H$ is
called the Hurst parameter, and may take any value in $(0,1)$.
Fractional Brownian motion $\fbm^H:(0,\infty) \to \mathbb{R}$ with
Hurst parameter $H$ satisfies
$$\mathbb{E}(\fbm_s^H \fbm_t^H) = \frac{1}{2} \big( \vert t \vert^{2H} +
\vert s \vert^{2H} - \vert t - s \vert^{2H} \big),$$
or, equivalently, 
$$\mathbb{E}(\vert \fbm_t^H - \fbm_s^H \vert^2) = \vert t-s\vert^{2H},
\,\,\,\,\, \fbm^H_0 = 0.$$ 
The process satisfies $\E \big( \fbm_t^H \big) = 0$ for all $t \in \R$.
When $H = 1/2$, it is Brownian motion, whose increments are of course independent. If $H
> 1/2$, increments over two given disjoint intervals are positively
correlated, while they have negative correlation if $H < 1/2$.

Fractional Brownian motion was first considered by Kolmogorov
\cite{kolmogorov} as a model of turbulence, and there is a now large
literature treating this family of processes, e.g. as models in
mathematical finance, and developing their stochastic calculus
\cite{booktwo}. The
stated characterization of fractional Brownian motion may permit
this family of processes to be considered as at least slightly
canonical. However, there are only a few examples of members of the
family arising as a scaling or continuum limit of a discrete model
(other than Brownian motion itself). The fractional Brownian motion
with Hurst parameter $1/4$ arises as a scaling limit for the tagged
particle in a one-dimensional symmetric exclusion process
\cite{arratia}. For general values of $H \in (0,1)$, fractional
Brownian motion has been exhibited as a scaling limit of an average
of a mixture of independent random walks, each walk having a decay
rate for the correlation of its increments, which rate is selected
independently from a law that depends on $H$ \cite{enriquez}. In
this paper, we present a simple discrete random walk that scales to
fractional Brownian motion. The process may be considered to be a
discrete counterpart to fractional Brownian motion.

In finance applications, as a model for the drift adjusted logarithm
of an asset price, fractional Brownian motion retains much of the
simplicity of ordinary Brownian motion (stationarity, continuity,
Gaussian increments) but dispenses with independence of increments,
thereby allowing for ``momentum effects'' (i.e., increment positive
correlations), which have been observed empirically in some markets.
One hypothesis is that momentum effects result from market
inefficiencies associated with insider trading; if an event occurs
of which the market is unaware, insider trading may cause the asset
price to change gradually, over a period of time, instead of all at
once --- initially because of trades by individuals with privileged
knowledge, later by a larger number of market participants.  These
ideas are discussed further in \cite{MR2087964}, where a model based
on shocks of this form is shown to be arbitrage free and and to have
fractional Brownian motion as a scaling limit. 
%(A strategy demonstrating arbitrage for fractional Brownian motion itself was given in \cite{rogers}, in an article that points out that, when fractional Brownian motion is represented as a convolution of white noise with a power-law kernel, it is the behaviour of this kernel near zero that is responsible for arbitrage; and this possibility is eliminated by smoothing that kernel near zero.)
(We will comment further on the possibility of arbitrage for fractional Brownian motion in Section \ref{secfurcom}.) See also
\cite{MR2281230}, which uses ``investor inertia'' to explain models
in which drift-adjusted logarithmic price is a stochastic integral
of fractional Brownian motion. Momentum effects appear naturally in
the random walks we introduce: i.e., it will be easy to see why an
``event'' (the sampling of an increment of the walk at one point in
time) has an influence on the expectation of future increments.

%We wonder whether the discrete process may be of interest to
%researchers in mathematical finance: if an event occurs of which the
%market is unaware, such as the breaking of an ice-cream production
%machine or a covert U.S. government decision to seek the overthrow
%of a Latin American government, insider trading may cause the price
%of a stock to jump at a sequence of times, initially because of
%trades by a few individuals with privileged knowledge, and then, by
%a larger number of market participants who hear rumors of the event.
%At least elements of this scenario may be modeled by a continuous
%time variant of the process we study, though no doubt it would be
%more natural to adapt the discrete process for a modeling
%application.

An informal description of the walk is as follows. Let $\mu$ be a
given law on the natural numbers (which we take to exclude zero).
The walk associated with law $\mu$ has increments, each of which is
either $-1$ or $1$. Independent samples $k_n$ of $\mu$ are attached
to the vertices of $\mathbb{Z}$. The sequence of increments $\big\{
\walk_k: k \in \mathbb{N} \big\}$ is such that, given the values $\big\{
\walk_k : k < n \big\}$, $\walk_n$ is set equal to the increment $\walk_{n -
k_n}$ obtained by looking back $k_n$ steps in the sequence. 
 The walk is then defined by adding the successive
increments from some fixed number. We
remark that the notion of determining the value $\walk_n$ of a process
at time $n$ by looking a random number of steps into the past also
appears in some of the urn models and reinforced random walk models
studied and surveyed by Pemantle in \cite{pemantle}.
Recall that in the traditional {\em P\'olya's urn}
process, the value of $\walk_n$ represents the color of the $n$th ball added to
an urn, and the conditional law of $\walk_n$ given the past is that of a uniform
sample from $\{\walk_1, \walk_2, \ldots, \walk_{n-1} \}$.  Our $\walk_n$ differs from the
P\'olya's urn process in that the past is infinite and the uniform sample is
replaced by one having a power law distance in time from the present.  

This description of the walk is a rough one, because it would require some
initial condition to construct it. We will turn to a Gibbs measure
formulation for a precise mathematical description.
Fixing the law $\mu$,
we write $G_{\mu}$ for the random directed spanning graph on
$\mathbb{Z}$ in which each vertex $z$ has a unique outward pointing
edge pointing to $z - k_z$, where $k_z$ is sampled independently
according to the measure $\mu$. We call $z - k_z$ the parent of $z$.
The ancestral line of $z$ is the decreasing sequence whose first
element is $z$ and each of whose terms is the parent of the previous
one.

A law $\lambda$ on functions mapping $\mathbb{Z}$ to $\{ -1,1 \}$ is
said to be a $\mu$-Gibbs measure provided that, for any
half-infinite interval $R = \big\{ x, \ldots, \infty \big\}$, the
conditional law of $\lambda$ on $R$, given its values in $R^c$, is
given by sampling $G_{\mu}$ and assigning to any element of $R$ the
value assigned to its most recent ancestor in $R^c$. We will write
$\walk = \walk_{\lambda}: \mathbb{N} \to \{ -1,1 \}$ for a realization of
the law $\lambda$. (In general, a
shift invariant probability measure on $\{-1,1\}^\mathbb{Z}$ that has a
specified conditional law for $X_0$ given $\{X_k: k < 0 \}$ is called a {\bf
$g$-measure} for that specification 
\cite{bramkal,kalikow,keane}.  One may view $\mu$-Gibbs
measures as $g$-measures for a particular specification.  We will not make
use of this more general framework here.)

A $\mu$-Gibbs measure is called {\bf
extremal} if it cannot be written as an average of two distinct
$\mu$-Gibbs measures.  There are at least two extremal $\mu$-Gibbs
measures, whatever the choice of $\mu$: those that assign unit mass
to either the constant function equal to $+1$ or to $-1$. If $\mu$
is such that the great common denominator of values in its support
is $1$, then the number of components in $G_{\mu}$ is either equal
almost surely to one or to $\infty$. (We furnish a proof of this
assertion in Lemma \ref{lemoneinf}.) If $G_{\mu}$ has infinitely
many components, then a one-parameter family $\big\{ \lambda_p : p
\in [0,1] \big\}$ of extremal Gibbs measures may be defined as follows.
To sample from $\lambda_p$, first sample $G_\mu$ and then
independently give each component of $G_\mu$ a value of $1$ (with
probability $p$) or $-1$ (with probability $1-p$).  Then assign all of the
vertices in that component the corresponding value.

Perhaps surprisingly, for such laws $\mu$, any
Gibbs measure is a mixture of these ones:
\begin{prop}\label{propone}
Let $\mu$ denote a probability measure on $\mathbb{N} = \big\{
1,2,\ldots \big\}$ the greatest common denominator of whose support is
equal to one.
If $G_{\mu}$ has one component almost surely, then
there are no extremal 
$\mu$-Gibbs measures other than the two trivial ones, $\lambda_0$
and $\lambda_1$. If $G_{\mu}$ has infinitely many components almost
surely, then the space of extremal 
Gibbs measures is the family $\big\{
\lambda_p : p \in [0,1] \big\}$.
\end{prop}
From this result, it is not hard to conclude that if the greatest common
denominator of the support of $\mu$ is some $k \not = 1$, and $\lambda$ is an
extremal $\mu$-Gibbs measure, then the restriction of $\lambda$ to points in $k
\mathbb Z+a$, for each $a \in \{ 0,1,\ldots, k-1 \}$, will be an extremal
$\mu$-Gibbs measure of the type described in Proposition 1, and that these
measures will be independent for different values of $a$ (though $p$ may depend
on $a$).  Hence, there is no real loss of generality in restricting to the
case that the greatest common denominator is $1$, as we do in Proposition 1
(and throughout most of the remainder of this paper).

Next, we define the measures $\mu$ that we will use for most of this
paper.
\begin{definition}\label{defone}
Let $\mu$ denote a probability measure on $\mathbb{N}$. For $\alpha
\in (0,\infty)$, we say that $\mu \in \G_{\alpha}$ if
there exists a slowly varying function $L:(0,\infty) \to
(0,\infty)$ for which
\begin{equation}\label{mudecay}
 \mu \big\{ n, \ldots,\infty \big\} = n^{- \alpha} L\big( n  \big)
\end{equation}
for each $n \in \mathbb{N}$.
Recall that by slowly varying is meant
\begin{equation}\label{slowvar}
\lim_{u \to \infty} \frac{L\big( u(1+r) \big)}{L(u)} = 1,
\end{equation}
for any $r > 0$.
\end{definition}

Note that if we required $L$ to be a constant function,
then the measures satisfying the first condition would be simply
those for which $\mu \{n,\ldots,\infty \}$ is a
constant times $n^{-\alpha}$.  The generalization to slowly varying
$L$ is quite natural, for the following reason. Let $R_\mu$
denote the random set $\big\{ \sum_{i=1}^j \walk_i : j \in \mathbb{N}
\big\}$ of values assumed by partial sums of an independent sequence
of samples $\walk_i$ of the law $\mu$.  (Clearly, $-1$ times the
ancestral line of $0$ has the same law as $R_\mu$.) Then it turns
out that when $\alpha \in (0,1)$, the random set $\epsilon R_\mu$
converges in law as $\epsilon \to 0$ to the range of a stable
subordinator with parameter $\alpha$ {\em if and only if} $\mu \in \G_{\alpha}$
(Theorem 8.3.1 of \cite{bgt}). Indeed, there is a sizable literature on probability
distributions with power law decays up to a slowly varying function
\cite{bgt}. However, throughout this paper, when a result is stated
for all $\mu \in \G_{\alpha}$, the reader may find it easier on a
first reading to focus on the special case that $L$ is constant. 

Our next result relates the decay rate of the tail of $\mu$
to the number of components in $G_{\mu}$.

\begin{prop}\label{proptwo}
Let $\mu \in \G_{\alpha}$ for some $\alpha \in (0,\infty)$. If
$\alpha > 1/2$, then $G_{\mu}$ almost surely has one component,
while, if $\alpha < 1/2$, then $G_{\mu}$ almost surely has
infinitely many.
\end{prop}

To any $\mu$-Gibbs measure $\lambda$, we associate a random walk $S
= S_{\lambda}: \mathbb{Z} \to \mathbb{Z}$ by setting $S(0) = 0$ and
$S(n) = \sum_{i=1}^n \walk
(i)$ for $n
> 0$, and $S(n) = - \sum_{i = -n}^{-1} \walk(i)$ for $n < 0$.
 We extend the domain of
$S$ to $\mathbb{R}$ by linearly interpolating its values between
successive integers.

We now ready to state our main result. Essentially, it says that, for $S$ the random walk associated to an extremal $\mu$-Gibbs measure for a choice of 
$\mu \in \Gamma_{\alpha}$, the time-scaled process $S(nt)$, further rescaled by subtracting its mean and multiplying by a deterministic $n$-dependent factor, converges to fractional Brownian motion with Hurst parameter $\alpha + 1/2$. This normalizing factor is written below as  $\tilde c
n^{-\frac12 - \alpha} L(n)$. The explicit form for the constant $\tilde c$ will be explained by the proof of Lemma \ref{lemvarsn}, in which the asymptotic variance of $S_n$ is determined.   
\begin{theorem}\label{thm}
For $\alpha \in (0,1/2)$, let $\mu \in \G_{\alpha}$. Let $L:(0,\infty) \to
(0,\infty)$ be given by (\ref{mudecay}).
Define $\varcon >0$ by means of
\begin{equation}\label{defvarcon}
\varcon^2 = \frac{\sum_{i=0}^{\infty} q_i^2}{2p(1-p)} \alpha \big( 2 \alpha + 1 \big) \Gamma \big( 1 - 2\alpha \big)^2 \Gamma \big( 2 \alpha \big) \cos \big( \pi \alpha \big),
\end{equation}
where $q_i = \mathbb{P} \big(i \in R_{\mu} \big)$ for $i \geq 1$, and $q_0 = 1$.
  Then, for each
$p \in (0,1)$, there exists a sequence of
couplings ${\rm C}_n$ of the process $$S_p^n: (0,\infty) \to
\mathbb{R}: t \to \varcon n^{-\frac{1}{2} - \alpha} L(n)  \Big(
S_{\lambda_p}\big( nt \big) - n(2p -1)t \Big)$$ and fractional
Brownian motion $\fbm_{\alpha + 1/2}$ with Hurst parameter $\alpha +
1/2$ such that, for each $T > 0$ and $\epsilon > 0$,
\begin{equation}\label{cncoupl}
 \lim_{n \to \infty} {\rm C}_n \Big(   \vert\vert  S_p^n - \fbm_{\alpha + 1/2} \vert\vert_{L_{\infty}\big( [0,T] \big)} > \epsilon  \Big)  = 0.
\end{equation}
\end{theorem}
\begin{subsection}{Discussion regarding Theorem \ref{thm}}
In light of Theorem 1.1, it is tempting to argue that the random walks we
construct are in some sense {\it the} canonical fractional analogs of simple
random walk on $\mathbb Z$.  To make this point, we note that any
random walk on $\mathbb Z$ (viewed as a graph) has increments in $\{-1,1\}$,
which means that its law, 
conditioned on $\{X_n: n < M \}$ for some fixed
$M$, is determined by the value of $E[\walk_n | \{\walk_k:
k<n\}]$.  If we assume stationarity of increments and we further posit that
$E[\walk_0 | \{\walk_k: k<0\}]$ has a simple form --- say, that it is a
monotonically increasing
linear function of $\{\walk_k: k<0\}$ --- then we have $$E[\walk_n | \{\walk_k: k<n\}] =
\sum_{i=1}^{\infty} p_i \walk_{n-i}$$ for some $p_i \geq 0$ with $\sum p_i \leq
1$.  It is not hard to show that if $\sum p_i < 1$, the process scales to
ordinary Brownian motion. We are therefore left with the case $\sum p_i =
1$, which corresponds to the walks we consider with $\mu(\{i\}) = p_i$.  

It seems plausible that any $\mu$ for which the conclusion of the theorem holds
(with the normalizing factor $\tilde c n^{-\frac12 - \alpha} L(n)$ replaced
by {\em some} deterministic function of $n$) must be a member of
$\Gamma_\alpha$.  Similarly, it seems highly plausible (in light of 
(\ref{nqtwo}))
that the closely related assertion that the variance of $S_n$ is
$n^{2\alpha+1}$ (multiplied by a slowly varying function) implies $\mu \in
\Gamma_\alpha$.  We will not prove either of these statements here.

We see that the model undergoes a phase transition at the value
$\alpha = \frac{1}{2}$. We remark that, while there are no
non-trivial $\mu$-Gibbs measures for $\mu \in \G_{\alpha}$ with
$\alpha > 1/2$, there is nonetheless a further phase transition at
the value $\alpha = 1$, which is the maximal value for which an
element $\mu \in \G_{\alpha}$ may have infinite mean. Indeed,
suppose that we define the walk associated to a measure $\mu$ by
instead specifying its domain to be the positive integers and then
making the following adjustment to the existing definition. As
previously, we take $\big\{ k_n : n \in \mathbb{N} \big\}$ to be a
sequence of independent samples of the measure  $\mu$, and, in the
case that $n - k_n \geq 0$, we continue to set the increment $\walk_n$
to be equal to $\walk_{n - k_n}$. In the other case, we choose this
increment to be $+1$ or $-1$ with equal probability, independently
of previous such choices. Then it is easily seen that the walk takes
infinitely many steps of each type if and only if $\mu$ has infinite
mean.

Finally, we mention that the decomposition of $\mathbb{Z}$ into components provided by
  Proposition \ref{proptwo} with a choice of $\mu \in \Gamma_\alpha$ with $\alpha
  < 1/2$ has something in common with the following process, discussed in
  \cite{aav}.  A system of particles, one at each element of $\Z$, are
  labelled, each by its location at an initial time $t = 0$. Each pair of adjacent
  particles consider swapping locations at an independent Poisson sequence
  of  times, but do so only if
  the higher-labelled particle lies on the right-hand-side in the
  pair just before the proposed swap. Each particle behaves as a second class particle in a totally
  asymmetric exclusion process, begun from an initial condition in which
  only sites to its left are occupied; as such, each has an asymptotic
  speed, which is a random variable having the uniform law on $[0,1]$. 
Defining a convoy
  to be the collection of locations inhabited at time zero by particles
  that share a common speed, (\cite{aav})
  shows that the convoy containing the origin is almost surely infinite,
  has zero density, and, conditionally on the speed of the particle
  beginning at the origin, is a renewal process. As such, the collection of
  convoys provides a partition of $\Z$ that has some similarities with
  $G_{\mu}$, when it has infinitely many components. 
\end{subsection}
\begin{subsection}{Further comments}\label{secfurcom}
 
In financial contexts, where fractional Brownian motion may serve as a model of evolving prices or preference strengths, it is particularly natural to aim to simulate its future trajectory given some fixed past. The discrete model that we introduce is very well suited to this problem, because, given the past, $N$ subsequent steps of the process may be sampled in time $O(N)$ (neglecting logarithmic corrections arising due to random number generation). In comparison, fast Fourier transform methods for approximately sampling fractional Brownian motion necessarily take at least $O(N \log N)$ time.    

Mandelbrot and Van Ness introduced a representation of fractional Brownian as a weighted average of white noise. The representation of fractional Brownian motion stated in footnote 3 on page 424 of \cite{mandvanness} inspires has a natural discrete counterpart, in which a weighted avarge of independent $\pm 1$ coin tosses is used.
However, the walks obtained in this way do not have $\pm 1$ increments and are not as simple to generate as those of this paper.
In regard to discretizations of fractional Brownian motion, we mention \cite{sottinen}, that introduces a random walk approximation and shows its weak convergence to fractional Brownian motion.

The question of whether (and what kind of) arbitrage is admitted in financial models based on fractional Brownian motion or geometric fractional Brownian motion has been explored at some length in the finance literature.  Fractional Brownian motion is not a semi-martingale: if the class of admissible trading strategies is large enough (and allows a potentially unbounded number of trades at arbitrarily small intervals) then one can indeed construct strategies that that admit positive return with probability one \cite{rogers}.  Even if such strategies involve trading that is too frequent to be practical, their existence complicates the construction of consistent derivative pricing models.  A good deal of literature has addressed ways in which transaction costs, other strategy restrictions, or minor modifications to the model can be imposed to make the arbitrage opportunites go away.  See \cite{bsv} for a survey of this literature that clarifies how the presence or absence of arbitrage depends on the precise choice of admissible trading strategies.
% for example \cite{MR2543527, MR2343807, MR2270939, MR2244029, MR2239592, MR2014249, MR1976518, MR2696663, MR1434408}.  
 We note that in a simple discrete model such as ours --- where the walk goes up or down by a unit increment at every discrete time step and all possible finite-length trajectories have positive probability --- one cannot have arbitrage in the strong sense, so these particular subtleties are not relevant to us.  It remains interesting, however, to think about the market efficiency implications of our model, and we present some preliminary thoughts on this subject in the appendix.
\end{subsection}
\begin{subsection}{Structure of the paper}
Section 2 is devoted to the analysis of $\mu$-Gibbs measures, and the
proofs of Propositions 1 and 2 are presented there. In section 3, we turn
to the convergence of the discrete process to fractional Brownian motion,
proving Proposition \ref{propfd}, which shows convergence in finite dimensional
distributions of the rescaled walk to the continuous process. 
A useful tool here is an explicit
asymptotic formula for the variance of the walk, which is presented in
Lemma \ref{lemvarsn}. In Section 4, we prove a correlation inequality
applicable to any extremal Gibbs measures in Proposition 4, and then apply
it to improve the topology of convergence to yield the
$L^{\infty}$-small coupling given in Theorem \ref{thm}.   

The discrete processes that we introduce form a natural counterpart to
fractional Brownian motion, and we are hopeful that they may serve as a
tool for the analysis of the continuous process. They also appear to have
various discrete relatives that are candidates for further enquiry. We take
the opportunity to discuss this by presenting several open problems, in Section 5.
We end with an appendix that describes potential economic and financial applications of the model.

\noindent{\bf Acknowledgments.}
We thank Yuval Peres for suggesting that the walks in question may converge to fractional Brownian motion. We thank S.R.S. Varadhan for discussions relating to regular variation.
\end{subsection}
\end{section}
\begin{section}{The space of Gibbs measures}
We now prove the first two propositions.
We begin with two simple lemmas.
\begin{lemma}\label{lemoneinf}
Let $\mu$ be a measure on the positive integers, the
greatest common denominator of its support being equal to one. Then
$G_{\mu}$ has either one component almost surely, or infinitely many
almost surely.
\end{lemma}
\noindent{\bf Proof.} We write $A_i$ for the ancestral line of $i \in \Z$.
Let $m_n = \mathbb{P} \big( A_0 \cap A_n \not= \emptyset \big)$. 
We now argue that if $\inf_n m_n
> 0$, then $G_{\mu}$ has one component almost surely.
If  $m : = \inf_n m_n > 0$, we may find $g_n \in \N$ such that the ancestral lines from
$0$ and $n$ meet at some integer exceeding $- g_n$ with probability at
least $m/2$. If the ancestral lines $A_0$ and $A_n$ fail to meet in
$\{-g_n,\ldots\}$, then we sample them until each first passes $-g_n$, at
$f_1$ and $f_2$, say. Without loss of generality, $f_1 < f_2$. 
Conditional on $A_0 \cap A_n \cap \big\{ -
g_n,\ldots \big\}= \emptyset$, and on $f_1$ and $f_2$, the probability that 
$A_0 \cap A_n \cap \big\{  f_1 - g_{f_2 - f_1},\ldots,f_1 \big\} \not= \emptyset$ is
at least $m/2$. Iterating, we construct a sequence of intervals each of
which has a probability at least $m/2$ of containing a point in $A_0 \cap A_n$,
conditionally on its precedessors not doing so. We
see that the ancestral lines indeed meet almost surely.

If
$m_n \to 0$ subsequentially, then we will show that
there exists a sequence $\big\{ n_i: i \in \mathbb{N} \big\}$ such
that
\begin{equation}\label{aiprob}
\mathbb{P} \Big( A_{n_i} \cap \bigcup_{j < i} A_{n_j} = \emptyset \Big)
 \geq 1 - i^{-2}.
\end{equation}
This suffices to show that $G_\mu$ has infinitely many components almost
surely, by the Borel-Cantelli lemma. To construct the sequence,
for $i \in \N$, set $q_i = \mathbb{P} \big( 0 \in A_i \big)$.
Then, whenever $n > l > 0$,
\begin{eqnarray}
& & \mathbb{P} \Big( A_0 \cap A_n \not= \emptyset \Big)
 \geq 
\mathbb{P} \Big( \Big\{ A_0 \cap A_n \not= \emptyset \Big\} \cap \Big\{   n - l \in
A_n  \Big\}  \Big) \nonumber \\
 & = &
\mathbb{P} \Big( A_0 \cap A_n \not= \emptyset  \Big\vert    n - l \in
A_n    \Big) \mathbb{P} \Big( n - l \in A_n \Big) =
\mathbb{P} \Big( A_0 \cap A_{n-l} \not= \emptyset  \Big) \mathbb{P} \Big( n
- l \in A_n \Big) \nonumber
\end{eqnarray}
so that 
\begin{equation}\label{mnqineq}
m_n \geq m_{n - l} q_l. 
\end{equation}
It follows readily from ${\rm g.c.d.} {\rm supp} (\mu) = 1$ that $q_j > 0$
for all sufficiently high $j \in \N$. 
Supposing that we have constructed 
an increasing sequence $\{ n_1,\ldots, n_{i-1} \}$ 
for which $q_{n_j} > 0$ for $j < i$, 
choose $n_i > n_{i-1}$ satisfying $q_{n_i} > 0$ and $m_{n_i} \sum_{j=1}^{i-1}
q_{n_j}^{-1} \leq i^{-2}$. Note then that
$$
 \mathbb{P} \Big( A_{n_i} \cap \bigcup_{j < i} A_{n_j} \not= \emptyset \Big)
\leq \sum_{j=1}^{i-1} \mathbb{P}  \Big( A_{n_i} \cap  A_{n_j} \not= \emptyset \Big)
 \leq m_{n_i} \sum_{j=1}^{i-1} q_{n_j}^{-1} \leq i^{-2},
$$
the second inequality by (\ref{mnqineq}). In this way, we construct a sequence
satisfying (\ref{aiprob}). $\Box$ \\
\noindent{\bf Proof of Proposition \ref{propone}.}
The case that $G_\mu$ has one
component almost surely is trivial.

Let $\mu$ be a probability measure on $\N$ whose support has greatest
common denominator one, and for which $G_{\mu}$ has infinitely many
components almost surely. Let $\lambda$ denote a $\mu$-Gibbs measure.
For $m \in \Z$, set $\sigma_m = \sigma \big\{ \walk_i: i < m \big\}$.
By the backwards martingale convergence theorem (Section XI.15 of \cite{doob}), 
\begin{equation}\label{pneqn}
p_{n,-\infty} = \lim_{m \to - \infty} \lambda \big( \walk_n = 1 \big\vert \sigma_m \big) \in \sigma_{-\infty}
 := \bigcap_{m < 0}\sigma_m
\end{equation}
exists $\lambda$-a.s. Note that it is an almost sure constant, because
$\lambda$ is extremal. 

We will now argue that 
\begin{equation}\label{pnmminusinf}
 p_{n,-\infty} = p_{m,-\infty} \qquad \textrm{$\lambda$-a.s.}
\end{equation}
We will denote the common value by $p$.

To show (\ref{pnmminusinf}), let $a,b \in {\rm supp}(\mu)$, $a \not=
b$. Consider firstly the case that $\vert b - a \vert$ divides $\vert n -m
\vert$. For $k \in \N$, we enumerate the ancestral line 
$A_k = \big\{ k = x_0(k),x_1(k),\ldots \big\}$ emanating from $k$ in decreasing order.
We define a coupling $\Theta$ of $A_n$ and $A_m$. At any given step $l \in
\N$, the initial sequences $\big\{ x_0(n),\ldots, x_l(n) \big\}$
and   $\big\{ x_0(m),\ldots, x_l(m) \big\}$ have been formed. If 
$x_l(m) \not= x_l(n)$, then we set $x_{l+1}(m) - x_l(m) = -q$, where $q$
has law $\mu$. If $q \not\in \{ a,b \}$, then we take  $x_{l+1}(n) - x_l(n) =
-q$ also. If $q \in \{ a,b \}$, then we take  $x_{l+1}(n) - x_l(n)$ to have
the conditional distribution of a sample of $\mu$ given that its value
belongs to $\{ a,b \}$. If $x_l(m) = x_l(n)$, we simply take $x_{l+1}(m) -
x_l(m) = x_{l+1}(n) - x_l(n) = -q$.

The sequence of differences $\big\{ x_l(m) - x_l(n): l \in \N \big\}$
performs a random walk with independent jumps, a jump equal to zero with
probability $1 - 2 \mu\big\{a \big\} \mu \big\{ b \big\}$, and, otherwise, with equal probability to
be $-\vert b - a \vert$ or $\vert b - a \vert$, until the walk reaches $0$.  The walk beginning at a
multiple of $\vert b - a \vert$, the sequence of differences reaches zero,
and then remains there, almost surely, by the recurrence of one-dimensional
simple random walk. 

If   $\vert b - a \vert$ does not divide $\vert n -m
\vert$, we begin by sampling initial segments of the two ancestral lines 
 $\big\{ x_0(n),\ldots, x_j(n) \big\}$
and   $\big\{ x_0(m),\ldots, x_j(m) \big\}$, where $j$ is the stopping time 
$$
j = \inf \Big\{ k \in \N: \vert b - a \vert \, \textrm{divides} \, 
 \vert x_k(n) -x_k(m) \vert  \Big\},
$$ 
which is almost surely finite, because ${\rm g.c.d.} {\rm supp} (\mu)
=1$. Indeed, this finiteness easily follows from the positivity of $q_j$ for all
sufficiently high $j$.
The coupling demonstrates that, given $\epsilon > 0$, there exists $l$
sufficiently negative that, $\sigma_l$-a.s.,
$$
 \Big\vert \lambda \Big( \walk_n  = 1 \Big\vert \sigma_l \Big)  - \lambda \Big( \walk_m = 1 \Big\vert
 \sigma_l  \Big) \Big\vert < \epsilon,
$$
from which (\ref{pnmminusinf}) follows.
 
We learn then that, for any $n \in \Z$,
$$
\lambda \big( \walk_n = 1  \big) = \E_{\lambda} \Pro \big( \walk_n = 1 \big\vert
\sigma_{-\infty} \big) = \E p_{n,-\infty} = p. 
$$

Let $L \subseteq \N$ be finite. Set 
$$
s_L = \lambda \Big( A_i, i \in L, \, \textrm{are not pairwise disjoint} \Big).
$$

For $m,l \in \Z$, $m < l$,
let $v_{l,m} = \sup A_l \cap \big\{ \ldots, m \big\}$ denote the first
element in the ancestral line of $l$ that is at most $m$.
For $L \subseteq \N$ finite with $m < \inf L$, write
$v_{L,m} = \big\{ v_{l,m}: l \in L \big\}$ for the locations of ancestors
of elements of $L$ that are at most $m$ whose child is at least $m+1$. 

For $\epsilon > 0$, set
$$
g_{\epsilon,L} = \sup \Big\{ m < \inf L : s_{v_{L,m}} < \epsilon \Big\}.
$$
Noting that $\big\{ \big\vert v_{L,m} \big\vert: m < \inf L \big\}$
is a non-decreasing sequence that assumes a common value for all sufficiently
high negative choices of $m$, we find that $g_{\epsilon,L} > - \infty$,
$\lambda$-a.s.

For $\epsilon > 0$ and $j,k \in \Z$, $k < j$, write
$$
\Omega_{k,j}^{(\epsilon)} = \Big\{
 \omega \in \{ -1,1 \}^{\{-\infty,\ldots,k \}}: 
 \lambda \Big( \walk_j = 1 \Big\vert \big( \ldots, \walk_k \big) = \omega \Big) 
\in \big( p - \epsilon , p + \epsilon \big) \Big\}.
$$
Note that, for all $\epsilon > 0$ and $j \in \N$,
\begin{equation}\label{laomega}
 \lim_{k \to - \infty} \lambda \Big( \Omega_{k,j}^{(\epsilon)} \Big) = 1,
\end{equation}
by (\ref{pneqn}) and (\ref{pnmminusinf}).
Given $L \subseteq \N$ finite, set
$$
j_{\epsilon,L} = \sup \Big\{ k < \inf L : 
  \lambda \Big( \Omega_{k,l}^{(\epsilon)} \Big) > 1 - \epsilon \, \,
 \textrm{for each $l \in L$} \Big\},
$$
so that $j_{\epsilon,L} > - \infty$ for each $\epsilon > 0$ and finite $L
\subseteq \N$, by (\ref{laomega}).

Let $k \in \Z$ and $\omega \in \{ -1,1 \}^{\{-\infty,\ldots,k\}}$.
Define $\big\{ Y_{k,l}^{(\omega)}: l \geq k + 1 \big\}$
by constructing an independent collection of ancestral lines (that do not
coalesce on meeting)
from the elements of $\big\{ l \in \Z: l \geq k + 1 \big\}$, stopping each
line on its arrival in $\big\{ l \in \Z: l \leq k \big\}$.
For each $l \geq k + 1$, we set $Y_{k,l}^{(\omega)}$ equal to the
$\omega$-value of the vertex at the end of the stopped ancestral line 
emanating from $l$.

Let $q_{\epsilon,L} = j_{\epsilon,v_{L,g_{\epsilon,L}}}$. Set
$\Omega_{\epsilon}^* = \bigcap_{l \in v_{L,g_{\epsilon,L}}}
\Omega_{q_{\epsilon,L},l}^{(\epsilon)}$, and note that
$\lambda \big( \Omega_\epsilon^* \big) \geq 1 - \vert L \vert \epsilon$.

Note further that, conditionally on $v_{L,g_{\epsilon,L}}$ and
for any $\omega \in  \{-1, 1 \}^{\{ - \infty, \ldots,
  q_{\epsilon,L} \}}$
for which 
$$
\lambda \Big( \walk_l = 1 \big\vert \big( \ldots,
\walk_{q_{\epsilon,L}} \big) = \omega \Big) \in \big( p -
\epsilon, p + \epsilon \big)
$$
for each $l \in v_{L,g_{\epsilon,L}}$, the set of values 
$\big\{ Y_{q_{\epsilon,L},l}: l \in v_{L,g_{\epsilon,L}}
\big\}$
are independent, each being equal to $+1$ with
probability at least $p - \epsilon$ and at most $p + \epsilon$. 
 
The labelling $\big\{ \walk_i: i \in L \big\}$ is determined by the set of
values $\big\{ \walk_i: i \in v_{L,g_{\epsilon,L}} \big\}$. However,
for each $\omega \in \{ -1,1 \}^{\{ \ldots,q_{\epsilon,L} \}}$
and for any finite set $\sigma \subseteq \Z$ for which
$v_{L,g_{\epsilon,L}} = \sigma$ is possible,
$$
{\rm TV}  \bigg( \Big(  \big( \walk_l: l \in v_{L,g_{\epsilon,L}} \big) \Big\vert
     v_{L,g_{\epsilon,L}} = \sigma, \big\{\ldots,\walk_{q_{\epsilon,L}} \big\}
  = \omega \Big), \big( Y_{j_{\epsilon,\sigma},l}^{\omega} : l \in
  \sigma \big) \bigg) \leq s_{\sigma}.
$$
Note that the right-hand-side 
is at most $\epsilon$, because $v_{L,g_{\epsilon,L}} = \sigma$ is possible.
By taking $\epsilon > 0$ arbitrarily small, we obtain the result.
$\Box$\\
\begin{definition}\label{defrmu}
Let $\mu$ denote a probability distribution on $\mathbb{N}$. Let
$R_{\mu}$ denote the random set $\big\{ \sum_{i=1}^j \walk_i : j \in
\mathbb{N} \big\}$ of values assumed by partial sums of an
independent sequence of samples $\walk_i$ of the law $\mu$. When such a
measure $\mu$ has been specified, we will write $p_n = \mu \{ n \}$
and $q_n = \mathbb{P} \big( n \in R_{\mu} \big)$. We adopt the convention
that $q_0 = 1$.
\end{definition}
\begin{lemma}\label{lemprobr}
Let $\alpha \in (0,1)$ and let $\mu \in \G_{\alpha}$.
\begin{enumerate}
\item
We have that
$$
 \sum_{n=0}^{\lfloor x \rfloor} q_n \sim \frac{1 - \alpha}{\Gamma(2 -
   \alpha) \Gamma(1 + \alpha)} x^{\alpha} L(x)^{-1},
$$
where $f \sim g$ means $\lim_{x \to \infty}{\frac{f(x)}{g(x)}} = 1$.
\item
Moreover, the sum
$$
 \sum_{n=0}^{\infty} q_n^2
$$
converges if $\alpha \in (0,1/2)$ and diverges if $\alpha \in
(1/2,1)$.
\end{enumerate}
\end{lemma}
\noindent{\bf Proof.}
Let $p,q:[0,\infty) \to [0,\infty]$, given by
$$
p(\lambda) = \sum_{n=0}^{\infty} p_n \exp \big\{ -\lambda n\big\}
$$
and
$$
q(\lambda) = \sum_{n=0}^{\infty} q_n \exp \big\{ -\lambda n\big\},
$$
denote the Laplace transforms of these two sequences. We will
analyze these transforms, (and, in later arguments, the Fourier
transforms) of such sequences by means of some Tauberian theorems.
A Tauberian theorem
asserts that a sequence (in our case) has a given rate of decay if
and only if its Laplace or Fourier transform has a corresponding
asymptotic behaviour close to zero. Additional hypotheses, such as
regular variation of the sequence (and, correspondingly, the
transform), are required.

We will make use of two basic results concerning slowly varying
functions. Firstly, the uniform convergence theorem (Theorem $1.2.1$ of 
  \cite{bgt})
states that, if $L:(0,\infty) \to (0,\infty)$ is slowly varying, then
\begin{equation}\label{uct}
\frac{L \big(  \lambda x \big)}{L \big( x \big)} \to 1,
\end{equation}
uniformly on each compact $\lambda$-set in $(0,\infty)$. Secondly, 
Potter's theorem (Theorem $1.5.6$ of \cite{bgt}) states in part that, for any slowly varying function 
 $L:(0,\infty) \to (0,\infty)$ that is bounded above and
below on any compact interval in $(0,\infty)$, and for any $\epsilon > 0$,
there exists $C =   C_\epsilon$ such that, for all $x,y > 0$,
\begin{equation}\label{lureps}
 \frac{L(y)}{L(x)} \leq C \max \Big\{ \big( y/x \big)^\epsilon , \big( y/x \big)^{-\epsilon} \Big\}. 
\end{equation}
From  $\big\{ p_n: n \in \mathbb{N} \big\} \in l_1$, the function $p$ is differentiable on $(0,\infty)$, with
\begin{equation}\label{dlamb}
\frac{d}{d \lambda} p(\lambda) = - \sum_{n=0}^{\infty} n p_n \exp \big\{ -\lambda n\big\}.
\end{equation}
Setting $U:(0,\infty) \to (0,\infty)$ according to $U(x) =
\sum_{n=0}^{\lfloor x \rfloor}n p_n$, it follows from 
$U(x) = 
\sum_{j=1}^{\lfloor x \rfloor} \sum_{i=j}^{\lfloor x \rfloor} p_i$,
(\ref{mudecay}),
 (\ref{uct}) and 
(\ref{lureps}) that
$U(x) \sim \frac{\alpha}{1-\alpha} x^{1-\alpha}L(x)$, $x \to \infty$.
A special case of Karamata's Tauberian theorem  (Theorem 1.7.1 of
\cite{bgt}) states that, if $\big\{ a_n : n \in \N \big\}$ is a sequence of
non-negative numbers, $l:(0,\infty) \to (0,\infty)$ is slowly varying,
and $c, \rho \geq 0$, then the following are equivalent:
\begin{eqnarray}
  \sum_{n=1}^{\lfloor x \rfloor}a_n   & \sim & \frac{c}{\Gamma(1+ \rho)}x^\rho
  l(x), \qquad x \to \infty, \nonumber \\
   \sum_{n=1}^\infty a_n e^{-\lambda n}    & \sim &  c \lambda^{-\rho}
  l\big ( \lambda^{-1} \big), \qquad \lambda \to 0^+. \nonumber
\end{eqnarray} 
Applying the theorem in the present case, we learn that
$$
   \sum_{n=1}^\infty n p_n e^{-\lambda n} 
  \sim  \frac{\alpha}{1 - \alpha} \Gamma \big( 2 - \alpha \big) 
 \lambda^{\alpha - 1} L\big( 1/\lambda \big), \qquad \lambda \to 0^+. 
$$
By (\ref{dlamb}), 
$1 - p(\lambda)  \sim \frac{\Gamma(2-\alpha)}{1 - \alpha} \lambda^\alpha
L\big( 1/\lambda \big)$, $\lambda \to 0^+$, since $\int_0^\lambda t^{\alpha
- 1} L \big( t^{-1} \big) dt \sim \alpha^{-1} \lambda^\alpha L \big(
\lambda^{-1} \big)$, by means of (\ref{uct}) and (\ref{lureps}). 
Noting that the Laplace transforms of the two sequences are related by $q = \big( 1 - p \big)^{-1}$, we learn that
\begin{equation}\label{qlaeqn}
q(\lambda)  \sim \frac{1 - \alpha}{\Gamma\big(2 - \alpha \big)} \lambda^{-
  \alpha} L\big( 1/\lambda \big)^{-1}, \qquad \lambda \to 0^+.
\end{equation}
The reverse implication in Karamata's Tauberian theorem then yields
the first statement of the lemma.

From (\ref{qlaeqn}) and
$$
q(\lambda) \leq \Big( \sum_{n=0}^{\infty} \exp \big\{ - 2 \lambda n \big\} \Big)^{1/2} \Big( \sum_{n = 0}^{\infty}{q_n^2} \Big)^{1/2} \leq \lambda^{-1/2} \Big( \sum_{n = 0}^\infty{q_n^2} \Big)^{1/2},
$$
(the first inequality is Cauchy-Schwarz; the second inequality holds
for $\lambda
> 0$ small enough), we find that $\sum_{n=0}^{\infty} q_n^2$
diverges, if $\alpha > 1/2$.

Set
$$
P(t) = \sum_{n=0}^{\infty} p_n \exp \big\{ itn \big\}
$$
and
$$
Q(t) = \sum_{n=0}^{\infty} q_n \exp \big\{ itn \big\}
$$
to be the Fourier transforms of the two sequences.

From $\mu \in \Gamma_\alpha$ with $\alpha \in (0,1)$, we learn from Theorem
1 of \cite{gelukdehaan} that 
$$
1 - {\rm Re} P(t) \sim \mu(t,\infty) \Gamma(1-\alpha) \cos \big( \alpha
\pi/2 \big) \qquad \textrm{as $t \downarrow 0$}
$$ 
and 
$$
{\rm Im} P (t) \sim \tan \big( \alpha \pi/2 \big) 
                   \Big( 1 - {\rm Re} P(t) \Big) \qquad \textrm{as $t \downarrow 0$,}
$$
whence 
\begin{equation}\label{pdecay}
 \Big\vert 1 - P(t) \Big\vert \sim \Gamma \big( 1 - \alpha \big) 
   t^{\alpha} L\big( t^{-1} \big) \qquad \textrm{as $t \downarrow 0$,}
\end{equation}
if $\mu (t,\infty) = t^{-\alpha} L(t)$. Similarly,
$$
 \Big\vert 1 - P(2\pi - t) \Big\vert \sim \Gamma \big( 1 - \alpha \big) 
   t^{\alpha} L\big( t^{-1} \big) \qquad \textrm{as $t \downarrow 0$.}
$$
From $Q = (1-P)^{-1}$, we find that, for any constant $C >
\Gamma(1-\alpha)^{-1}$, there exists $c > 0$ such that
\begin{equation}\label{gtalbd}
 \vert Q(t) \vert \leq C t^{-\alpha} L(1/t)^{-1}
\end{equation}
for $0 < t < c$, and, similarly,
\begin{equation}\label{twopiq}
 \vert Q\big(2\pi -  t \big) \vert \leq C t^{-\alpha} L(1/t)^{-1}.
\end{equation}
From the support of $\mu$ having greatest common denominator one,
$\vert P(\lambda) \vert < 1$ for $\lambda \in (0,2\pi)$. By the
continuity of $P:(0,\infty) \to \mathbb{C}$,
\begin{equation}\label{pcpt}
\sup \Big\{  \vert P(\lambda) \vert  : \lambda \in [c,2\pi - c] \Big\}
 \leq 1 - c,
\end{equation}
for $c > 0$ small.

By Parseval's identity,
$$
 \sum_{n=0}^{\infty} q_n^2 = \frac{1}{2\pi}\int_0^{2\pi} \vert Q(t) \vert^2 dt,
$$
which, by (\ref{gtalbd}), (\ref{twopiq}) and (\ref{pcpt}) is finite,
provided that $\alpha < 1/2$. $\Box$ \\
\noindent{\bf Proof of Proposition \ref{proptwo}.}
Let $\alpha > 0$ and $\mu \in \G_{\alpha}$.
Let $n,m \in \mathbb{N}$ satisfy $n < m$. We wish to show that,
if $\alpha > 1/2$, then
 $n$ and $m$ lie in the same component of $G_{\mu}$,
while, in the case that $\alpha < 1/2$,
 there is a positive probability (depending on $m-n$) that they lie in different components.

Let $R^{(1)}_\mu$ and $R^{(2)}_\mu$ be independent samples of $R_{\mu}$. It is easily seen that $n$ and $m$ almost surely lie in the same component of $G_{\mu}$ if and only if $\mathbb{E} \big( \big\vert \big\{ n - R^{(1)}_{\mu} \big\} \cap \big\{ m - R^{(2)}_{\mu} \big\} \big\vert \big) = \infty$.
Moreover, if there exists some pair $(m,n) \in \mathbb{Z}^2$ for which this expectation is infinite, then clearly, it is infinite for all such pairs. Hence, there exists one component of $G_{\mu}$ precisely when $\mathbb{E} \big\vert R_{\mu}^1 \cap R_{\mu}^2 \big\vert$ is infinite. This expectation is
equal to $\sum_{n=1}^{\infty} q_n^2$, so the result follows from
Lemma \ref{lemprobr}(ii). $\Box$
\end{section}
\begin{section}{Convergence to fractional Brownian motion}
In this section, we establish that the walk has a fractional Brownian motion scaling limit in the sense of finite dimensional distributions:
\begin{prop}\label{propfd}
Let $\mu \in \G_{\alpha}$, 
for some $\alpha \in (0,1/2)$. Recall the constant $\varcon = \varcon(p)$ for $p \in (0,1)$ from (\ref{defvarcon}). Then, for each $p \in (0,1)$,
the process $(0,\infty) \to \mathbb{R}: t \to c n^{-\frac{1}{2} - \alpha} L(n)
\Big( S_{\lambda_p}\big( nt \big) - n(2p -1)t \Big)$
converges weakly (in the sense of finite dimensional distributions)
as $n \to \infty$ to fractional Brownian motion with Hurst parameter $\alpha + 1/2$.
\end{prop}
In preparation for the proof,
let $X: \mathbb{Z} \to \{ -1,1 \}$ be a sample of $\lambda_p$, and let $S:
\mathbb{Z} \to \mathbb{Z}$ be the random walk satisfying $S(0) = 0$ and $S_n - S_{n-1} = \walk_n$ for each $n \in \mathbb{Z}$.
It is our aim to show that $S_n$ is approximately normally distributed, when $n$ is chosen to be high.
We begin by finding an explicit expression for the variance of $S_n$.
\begin{lemma}\label{lemvarsn}
We have that
$$
 {\rm Var} \big( S_n \big) \sim
 \frac{4p(1-p)K_\alpha}{\big( \vert Q \vert^2 \big)_0} n^{2\alpha + 1} L(n)^{-2},
$$
where
$\big( \vert Q \vert^2 \big)_i = \sum_{j=0}^\infty q_j q_{i+j}$
denotes the $i$-th Fourier coefficient of $\vert Q \vert^2$, and where
$$K_\alpha = 
 \frac{1}{2\alpha(2\alpha+1)} \Big( \Gamma \big( 1-2\alpha \big)^2 \Gamma\big( 2\alpha \big)
   \cos\big( \pi \alpha \big) \Big)^{-1}. 
$$
\end{lemma}
\noindent{\bf Proof.}
We begin by showing that
\begin{equation}
 {\rm Var} S_n = \frac{4p(1-p)}{\big( \vert Q \vert^2 \big)_0}
  \bigg( 2 \sum_{i=1}^n \big( n - i  \big) \big( \vert Q \vert^2
  \big)_i  \, +  \, n \big( \vert Q \vert^2 \big)_0    \bigg). \label{varsncom}
\end{equation}
(In fact, (\ref{varsncom}) holds, more
generally, for any $\mu$ for which $G_\mu$ has a.s.\ infinitely many
components.)
To do so, write $T_1,\ldots T_{r(n)}$ for the trees having non-empty
intersection with $\big\{ 1,\ldots, n \big\}$. Then
\begin{eqnarray}
{\rm Var} \big( S_n  \big) & = & 4p(1-p) \E \sum_{i=1}^{r(n)} 
 \Big\vert T_i \cap \big\{ 1 ,\ldots, n\big\} \Big\vert^2 \nonumber \\
 & = & 4p(1-p) \sum_{i=1}^n \sum_{j=1}^n \Pro \Big( A_i \cap A_j \not=
 \emptyset \Big). \nonumber
\end{eqnarray}
Note that, for $i < j$,
$$
\E \Big( A_i \cap A_j \Big) = \sum_{k = 0}^\infty q_k q_{j-i+k} = \big(
\vert Q \vert^2  \big)_{j-i},
$$
and also that
$$
\E \big( A_i \cap A_j \big) = \Pro \big( A_i \cap A_j \not= \emptyset \big)
 \sum_{i=0}^\infty q_i^2  =  \Pro \big( A_i \cap A_j \not=
 \emptyset \big) \big( \vert Q \vert^2 \big)_0 ,
$$
since we adopt the convention that $q_0 = 1$.
Hence,
$$
{\rm Var} \big( S_n  \big)   =   
 \frac{4p(1-p)}{\big( \vert Q \vert^2 \big)_0}
  \sum_{i=1}^n \sum_{j=1}^n 
 \big( \vert Q \vert^2 \big)_{\vert j - i \vert},
$$
whence, (\ref{varsncom}). By (\ref{varsncom}), it suffices to show that
\begin{equation}\label{nqtwo}
\sum_{i=1}^n \big( n  - i \big)
 \big( \vert Q \vert^2 \big)_i \sim
 K_\alpha n^{2\alpha + 1} L(n)^{-2},
\end{equation}
because Lemma \ref{lemprobr}(ii) implies that
$\big( \vert Q \vert^2 \big)_0 < \infty$. To this end, note that, 
by (\ref{pdecay}) and $Q = (1 - P)^{-1}$,
$$
\big\vert Q(t) \big\vert \sim \Gamma(1-\alpha)^{-1} t^{-\alpha} L\big( t^{-1} \big)^{-1},
$$
so that
$$
 \vert Q(t) \vert^2 \sim  \Gamma(1-\alpha)^{-2} t^{-2\alpha} L\big( t^{-1} \big)^{-2},
$$
By Theorem 4.10.1(a) of \cite{bgt}, using $\vert Q(t) \vert^2 \in \R$,
$$
\sum_{i=1}^n \big( \vert Q(t) \vert^2 \big)_i
   \sim  \Big( \Gamma \big( 1-2\alpha \big)^2 \Gamma\big( 2\alpha \big)
   \cos\big( \pi \alpha \big) \Big)^{-1}
 \frac{n^{2\alpha} L(n)^{-2}}{2\alpha}.
$$
It follows that
$$
 \sum_{i=1}^n \big( n  - i \big)
  \big( \vert Q \vert^2 \big)_i 
  \sim \Big( \Gamma \big( 1-2\alpha \big)^2 \Gamma\big( 2\alpha \big)
   \cos\big( \pi \alpha \big) \Big)^{-1} 
n^{2\alpha + 1}L(n)^{-2} \int_0^1 (1 -  x) x^{2\alpha - 1} dx.
$$
Evaluating the integral, we obtain (\ref{nqtwo}). $\Box$

%To sketch the proof, we may consider the sum $S_n$ as being
%determined successively: suppose that we know the values of $\walk_i$
%for $i \leq M$, with $M$ large and negative. We may write an
%estimate for $S_n$ on the basis of this information, using the mean
%value of $S_n$ given these values of $\walk_i$. We successively increase
%the value of $M$ towards zero, refining our estimate of $S_n$. The
%updating estimate is a martingale, and the variance of each revision
%in the estimate is much smaller than the variance of the total. Our
%eventual estimate is correct, and the martingale central limit
%theorem applies to show that $S_n$ has an approximately normal
%distribution.

\noindent{\bf Proof of Proposition \ref{propfd}.}
There being a unique stationary Gaussian process $\fbm^H:\mathbb{R} \to \mathbb{R}$ with covariances given by $\E \vert \fbm_t^H \vert^2 = t^{2\alpha + 1}$, it suffices, in light of Lemma \ref{lemvarsn}, to establish that $S_n$ has a distribution that is asymyptotically Gaussian.  
We will make use of the following notation.
For $M \in \mathbb Z$, set $\sigma_M = \big\{ \walk_i:i< M \}$.
Let $\walk_n^M$ and $S_n^M$ denote the expected
values of $\walk_n$ and $S_n$ respectively given $\sigma_M$.  (Note that $\walk_n^M$ is the expected value of
$\walk_k$ where $k$ is the first element in the ancestral line of $n$
that is less than $M$.) Clearly, $\walk_n^M$ and $S_n^M$ are martingales
in $M$, and $S_n^M = S_n$ when $n \geq M$. We will establish that $S_n$ has an asymptotically Gaussian law by applying the martingale central limit theorem to
$S_n^M$. This will require showing that $S_n^M$ has small increments
compared to its total size and a use of a long-range
near-independence argument to show
that the sum of the conditional variances of the increments is
concentrated. We now state the martingale central
limit theorem in the form that we require. See Theorem 7.2 in Chapter 7 of
\cite{durr}, and the remark following its proof, for a derivation.
\begin{definition}
We say that $\walk_{n,m}, \F_{n,m}$, $n \geq 1$, $1 \leq m \leq k_n$
is a martingale difference array if
$\walk_{n,m} \in \F_{n,m}$ and $\E \big( \walk_{n,m} \big\vert \F_{n,m-1} \big) =
0$
for $1 \leq m \leq k_n$. Let $V_n = \sum_{i=1}^{k_n}\E \big( \walk_{n,m}^2
\big\vert \F_{n,m-1} \big)$.
\end{definition} 
\begin{theorem}\label{martclt}
Suppose that $\big\{ \walk_{n,m} , \F_{n,m}  \big\}$ is a martingale difference
array. Let $\walk_{n,0} \in \F_{n,0}$, $n \in \N$. Set $S_n = \sum_{j=0}^{k_n}
\walk_{n,j}$. Assume that, for some sequence $\epsilon_n \to 0$,
\begin{itemize}
\item  $\big( \E V_n \big)^{-1} \big\vert \walk_{n,m} \big\vert^2 \leq \epsilon_n \, \,
  \textrm{for all $m \in \{1,\ldots, k_n  \}$}$
\item $\big( \E V_n \big)^{-1} {\rm Var} \big( \walk_{n,0} \big) \leq
  \epsilon_n$ for $n$ sufficiently high,
\item $\frac{V_n}{\E V_n} \to 1, \, \, \textrm{in probability.}$ 
\end{itemize}
Then $\E(V_n)^{-1/2} S_n$ converges in distribution to the normal
distribution of mean zero and unit variance.
\end{theorem} 
We will apply the result in the following way. 
Let $\big\{ \epsilon_n: n \in \N \big\}$ be a given sequence converging to
zero. Let $\big\{ k_n: n \in \N \big\}$ be a sequence satisfying 
\begin{equation}\label{varsnkn}
{\rm Var} \big( S_n^{-k_n} \big) \leq 2^{-1} \epsilon_n
        \frac{4p(1-p)K_\alpha}{(\vert Q \vert^2 )_0} n^{2\alpha + 1}
        L(n)^{-2}, 
\end{equation}
where the constant $K_\alpha$ is specified in Lemma \ref{lemvarsn}. 
We in addition assume that $k_n \geq Cn$ for each $n \in \N$, and  for some
large constant $C > 0$.
We 
choose as martingale difference array, $\walk_{n,0} = S_n^{-k_n}$, $\F_{n,0} = 
\sigma_{- k_n}$ for each $n \geq 1$, and 
$\walk_{n,i} = S_n^{-k_n + i} - S_n^{-k_n + i - 1}$, 
$\F_{n,i} = \sigma_{- k_n + i}$ for each $n \geq 1$
and $1 \leq i \leq k_n + n$.

Clearly, we must verify that the hypotheses of Theorem \ref{martclt}
hold. The first column of the martingale difference array has been chosen
to meet the second condition listed in the theorem, as we now confirm. 
Note that
\begin{equation}\label{sumofvar}
 \E V_n  = {\rm Var} \Big(  S_n - S_n^{k_n}  \Big).
\end{equation}
By (\ref{varsnkn}) and Lemma \ref{lemvarsn}, we see that 
\begin{equation}\label{varsnf}
 {\rm Var} \big( S_n^{- k_n } \big) \leq 
2^{-1} \epsilon_n {\rm Var} S_n \Big( 1 + o(1) \Big).
\end{equation}
Using ${\rm Var} (S_n) = {\rm Var}(S_n - S_n^{-k_n}) + 
{\rm Var} (S_n^{-k_n})$
and (\ref{sumofvar}), we obtain 
${\rm Var} S_n^{-k_n} \leq \epsilon_n \E V_n$ for $n$ sufficiently high,
which is indeed the second
hypothesis of Theorem \ref{martclt}.

The general martingale difference term takes the form $S_n^m -
S_n^{m-1}$, for some $n \geq 1$ and $m \leq n$. A basic observation is
that, for any $M < n$,
\begin{equation}\label{xupdate}
 \walk_n^{M+1} - \walk_n^M = \big( \walk_M - \walk_M^M \big) q_{n-M},
\end{equation}
where $q_i$ (as in Definition \ref{defrmu}) denotes the probability
that the vertex $i$ has $0$ as an ancestor.  As stated earlier,
$\walk_n^M$ is the expected value of the value observed by tracing back
the ancestral line in a sample of $G_{\mu}$ from $n$, until an
ancestor with index strictly less than $M$ is reached. The
difference in the above equation, then, is due to the event that the
ancestral line reaches $M$.
 
If $n \geq 1$ and $m < 0$, note that $S_n^{m+1} - S_n^m = \sum_{i=1}^n
\big( \walk_i^{m+1} - \walk_i^m \big)$, so that (\ref{xupdate}) yields 
\begin{equation}\label{formneg}
S_n^{m+1} - S_n^m =  \big( \walk_m - \walk_m^m \big) F_n^m, 
\end{equation}
where we define $F_n^m = \sum_{i=1}^n q_{i - m}$. 
In the case where 
$n \geq 1$ 
and 
$m \in \{0,\ldots,n-1\}$, 
on the other
hand, 
$$
S_n^{m+1} - S_n^m = \sum_{i=m}^{n-1} \big( \walk_i^{m+1} - \walk_i^m \big) =
\big( \walk_m - \walk_m^m \big) \sum_{i=m}^{n-1} q_{i - m},
$$ 
by the convention
that 
$q_0 = 1$. 
We record this as
\begin{equation}\label{formpos}
S_n^{m+1} - S_n^m = \big( \walk_m - \walk_m^m \big) \sum_{i=0}^{n-m-1} q_i .
\end{equation}
Note also that, for any $m \in \Z$,
$$
{\rm Var} \Big( \walk_m - \walk_m^m \Big\vert \sigma_m \Big) = 4 P_m \big( 1 - P_m \big),
$$
where $P_m$ is defined to be $\mathbb{P} \big( \walk_m = 1 \big\vert \sigma_m
\big)$ and may be written $\big( \walk_m^m + 1 \big)/2$.
The definition of $V_n$ and the expressions (\ref{formneg}) and
(\ref{formpos}) now give the formula
$$
V_n  =  
\sum_{M = - k_n}^{-1}   \big( F_M^n
\big)^2 {\rm Var} \big( \walk_M - \walk_M^M \big\vert \sigma_M \big)  +
\sum_{m=1}^n \Big( \sum_{l=1}^{n-m} q_l \Big)^2  
{\rm Var} \Big( \walk_m - \walk_m^m \Big\vert \sigma_M \Big), 
$$
which may be written
\begin{equation}\label{vnrep}
V_n  = 
4 \sum_{M = - k_n}^{-1}  P_M(1-P_M) \big( F_M^n
\big)^2  + 4 \sum_{m=1}^n P_m(1-P_m) \Big( \sum_{l=1}^{n-m} q_l \Big)^2. 
\end{equation}
As such, the following lemma shows that the third hypothesis of 
Theorem \ref{martclt} is satisfied in the present case.
\begin{lemma}\label{lemvn}
Setting $c_0 = 4 \E \big( P_M (1-P_m) \big)$,
$$
 4 \sum_{M = - k_n}^{-1}  P_M(1-P_M) \big( F_M^n
\big)^2 = c_0 \sum_{M = - k_n}^{-1}   \big( F_M^n
\big)^2 \, \Big( 1 + E_1(n) \Big),  
$$
where $E_1(n) \to 0$ in probability. We have further that
$$
\sum_{m=1}^{n} P_m (1 - P_m) \Big( \sum_{l=1}^{n-m}q_l\Big)^2  =  c
\sum_{m=1}^{n} \Big( \sum_{l=1}^{n-m}q_l\Big)^2 \, \Big( 1 + E_2(n) \Big) .
$$
where $E_2(n) \to 0$ in probability.
\end{lemma}
The ergodic theorem might be used to prove a concentration inequality of
this sort. We prefer, however, to derive it directly, from a second moment
estimate. We begin by observing the following.
\begin{lemma}\label{leml}
Set
$$
Z_i =  {\rm Var}\big( \walk_i - \walk_i^i \big\vert \sigma_i \big) =
4P_i(1-P_i),
$$
and write
$$
\rho_{l,m} = \mathbb{E} \big(  Z_l Z_m \big) - \mathbb{E} \big(  Z_l \big) \mathbb{E} \big( Z_m \big).
$$
Then, for each $\epsilon > 0$, there exists $K_0 \in \mathbb{N}$
such that $\vert l - m \vert \geq K_0$ implies that
$\vert \rho_{l,m} \vert < \epsilon$.
\end{lemma}
\noindent{\bf Proof.}
For $i \in \Z$ and $k \in \N$, set
$\walk_{i,k}^i = \sum_{j=1}^k \mu \{ j \} \walk_{i-j}$. Note that
\begin{equation}\label{xiineq}
\big\vert \walk_i^i - \walk_{i,k}^i \big\vert \leq
 \sum_{j = k+1}^\infty \mu \{ j  \} = \big( k + 1 \big)^{-\alpha} L (k+1).
\end{equation}
For $i,j \in \Z$ and $k \in \N$, set
$$
R_{i,j,k} = \Big\{  \bigcup_{m = i - k}^{i-1} A_m  \, \cap \, 
      \bigcup_{n = j - k}^{j-1} A_n = \emptyset   \Big\},
$$
where recall that we write $A_n$ for the ancestral line of $n \in \Z$.
Shortly, we will show that, for each $\epsilon > 0$ and $k \in \N$,
there exists $n_0 \in \N$ such that, if $i,j \in \Z$ satisfy
$\vert i - j \vert \geq n_0$, then
\begin{equation}\label{riineq}
 \lambda \Big( R_{i,j,k} \Big) < \epsilon.
\end{equation}
For now, we show that this suffices for the proof of the lemma. 

Note that, for any $i,j \in \Z$ and $k \in \N$ for which $\vert i -j \vert
\geq k$, if $Y \in 
\sigma \big\{ \walk_{i-k},\ldots, \walk_{i-1} \big\}$ and $Z \in 
\sigma \big\{ \walk_{j-k},\ldots, \walk_{j-1} \big\}$
\begin{equation}\label{xrcon}
 \E \Big( Y Z \Big\vert R_{i,j,k} \Big)
 =  \E \Big( Y  \Big\vert R_{i,j,k} \Big) 
\E \Big( Z \Big\vert R_{i,j,k} \Big).
\end{equation}
Note that
\begin{eqnarray}
 & & \E \Big( P_i(1 - P_i) P_j(1 - P_j)  \Big) \nonumber \\
 & = & \E   \Big( \frac{1 + \walk_i^i}{2}  \Big)
  \Big( \frac{1 - \walk_i^i}{2}  \Big)  \Big( \frac{1 + \walk_j^j}{2}  \Big)  \Big(
  \frac{1 - \walk_j^j}{2}  \Big) \nonumber \\ 
 & = & \E   \Big( \frac{1 + \walk_{i,k}^i}{2}  \Big)
  \Big( \frac{1 - \walk_{i,k}^i}{2}  \Big)  \Big( \frac{1 + \walk_{j,k}^j}{2}  \Big)  \Big(
  \frac{1 - \walk_{j,k}^j}{2}  \Big)  + O\Big(k^{- \alpha}L(k)\Big) \nonumber \\ 
 & = & \E \bigg(  \Big( \frac{1 + \walk_{i,k}^i}{2}  \Big)
  \Big( \frac{1 - \walk_{i,k}^i}{2}  \Big)  \Big( \frac{1 + \walk_{j,k}^j}{2}  \Big)  \Big(
  \frac{1 - \walk_{j,k}^j}{2}  \Big) \Big\vert R_{i,j,k} \bigg) + O(\epsilon) +
  O\Big(k^{- \alpha}L(k)\Big) \nonumber \\
& = & \E \bigg(  \Big( \frac{1 + \walk_{i,k}^i}{2}  \Big)
  \Big( \frac{1 - \walk_{i,k}^i}{2}  \Big) \Big\vert R_{i,j,k} \bigg)
 \E \bigg(  \Big( \frac{1 + \walk_{j,k}^j}{2}  \Big)  \Big(
  \frac{1 - \walk_{j,k}^j}{2}  \Big) \Big\vert R_{i,j,k} \bigg) + O(\epsilon) +
  O\Big(k^{- \alpha}L(k)\Big) \nonumber \\
& = & \bigg( \E \Big( \frac{1 + \walk_{i,k}^i}{2}  \Big)
  \Big( \frac{1 - \walk_{i,k}^i}{2}  \Big) + O\Big( \lambda\big( R_{i,j,k}^c
  \big) \Big) \bigg) \nonumber \\
 & & \qquad
  \bigg( \E  \Big( \frac{1 + \walk_{j,k}^j}{2}  \Big)  \Big(
  \frac{1 - \walk_{j,k}^j}{2}  \Big) + O\Big( \lambda\big( R_{i,j,k}^c
  \big) \Big) \bigg) + O(\epsilon) +
  O\Big(k^{- \alpha}L(k)\Big) \nonumber \\
 & = & \E \Big( \frac{1 + \walk_i^i}{2}  \Big)
  \Big( \frac{1 - \walk_i^i}{2}  \Big) 
\E  \Big( \frac{1 + \walk_j^j}{2}  \Big)  \Big(
  \frac{1 - \walk_j^j}{2}  \Big) + O(\epsilon) +
  O\Big(k^{- \alpha}L(k)\Big), \nonumber
\end{eqnarray}
the second equality by (\ref{xiineq}), the third by $\vert i - j \vert \geq
n_0$ and (\ref{riineq}), the fourth by (\ref{xrcon}), and the sixth by
(\ref{xiineq}) and (\ref{riineq}).
Thus, for any $\epsilon > 0$, there exists $n_0 \in \N$ such that, if $i,j
\in \Z$ satisfy $\vert i - j \vert \geq n_0$, then
$$
{\rm Cov} \Big( P_i(1-P_i) , P_j(1 - P_j) \Big) < \epsilon.
$$
This completes the proof of the lemma, subject to verifying (\ref{riineq}).
Note that, for this, it suffices to show that, for any $\epsilon > 0$,
there exists $n_0 \in \N$ such that, if $i,j \in \Z$ satisfy $\vert i - j
\vert \geq n_0$, then $\lambda \big( A_i \cap A_j \not= \emptyset \big) <
\epsilon$. To see this, note that, for $i < j$,
$$
\lambda \big( A_i \cap A_j \not= \emptyset \big) 
 \leq  \mathbb{E} \big( A_i \cap A_j \big) =
\sum_{n=0}^{\infty} q_n q_{n + i - j} \leq
 \Big( \sum_{n=0}^{\infty} q_n^2 \Big)^{1/2}  \Big( \sum_{n=0}^{\infty} q_{n + i - j}^2 \Big)^{1/2},
$$
so that Lemma \ref{lemprobr}(ii)
yields the desired conclusion. $\Box$ \\
\noindent{\bf Proof of Lemma \ref{lemvn}.}
To verify the first claim, we aim to show that
\begin{equation}\label{vyn}
 \textrm{Var} Y_n = o \Big(  \big( \mathbb{E} Y_n \big)^2 \Big),
\end{equation}
where here we write
$$
Y_n = 
 4 \sum_{M = - k_n}^{-1}  P_M(1-P_M) \big( F_M^n \big)^2.
$$
We will require 
\begin{equation}\label{fmn}
 \frac{\max_{- k_n <M<0} \big( F_M^n \big)^2}{ \sum_{M =
     - k_n}^{-1} \big( F_M^n \big)^2} \to 0
\end{equation}
as $n \to \infty$.

Indeed, by computing the second moment of $Y_n$, (\ref{vyn}) follows
from (\ref{fmn}) by means of Lemma \ref{leml}.

We turn to proving (\ref{fmn}). Recall that we 
supposed that $k_n \geq Cn$, with $C > 0$ a large constant. Note that
$$
\sum_{M =- k_n}^{-1} \big( F_M^n \big)^2 =
\sum_{m=1}^{k_n}
\Big( \sum_{i=1}^n q_{i+m} \Big)^2
 \geq \sum_{m=1}^n \Big( \sum_{i=1}^n q_{i+m} \Big)^2
 \geq c n^{2\alpha + 1} L(n)^{-2},
$$
with $c > 0$ a constant satisfying $c < \Big( \frac{2(1-\alpha)}{\Gamma(2 -
  \alpha) \Gamma(1 + \alpha)} \Big)^2$, and 
the latter inequality by Lemma \ref{lemprobr}(i),
(\ref{uct}) and 
(\ref{lureps}).
For $0 \leq m \leq \lfloor Cn \rfloor$ and a constant 
$C_0 >  \frac{2(1-\alpha)}{\Gamma(2 - \alpha) \Gamma(1 + \alpha)}$
$$
\sum_{i=1}^n q_{i+m} \leq C_0 n^{\alpha} L(n)^{-1},
$$
by  Lemma \ref{lemprobr}(i) and (\ref{uct}). We learn that
\begin{equation}\label{nmzerobd}
 \frac{\max_{- k_n < M< 0} \big( F_M^n \big)^2}{ \sum_{M = -\lfloor Cn \rfloor}^{-1} \big(
   F_M^n \big)^2} \leq  C^2 c^{-1} n^{-1}.
\end{equation}
To treat the bound on $F_M^n$ for $- k_n < M < - n$, note that, if $m > n$,
then either
$$
\sum_{i=1}^n q_{i+j} \geq 2^{-1} \sum_{i=1}^n q_{i+m} \, \, \, 
\textrm{for each $j \in \big\{ m - n/2, \ldots, m  \big\}$}
$$
or
$$
\sum_{i=1}^n q_{i+j} \geq 2^{-1} \sum_{i=1}^n q_{i+m} \, \, \, 
\textrm{for each $j \in \big\{ m + 1, \ldots, m + n/2 \big\}$},
$$
the first alternative holding provided that at least one half of the sum 
$\sum_{i=1}^n q_{i+m}$ is contained in its first $\lfloor n/2 \rfloor$
terms.
From this, we conclude that, for any $M < -n$,
$$
\frac{\big( F_M^n \big)^2}{\sum_{M =- \infty}^{-1} \big( F_M^n \big)^2}
\leq 8n^{-1} \Big( 1 + o(1) \Big),
$$
whence
$$
\frac{\big( F_M^n \big)^2}{\sum_{M =- k_n}^{-1} \big( F_M^n \big)^2}
\leq 10 n^{-1},
$$
by increasing the value of $k_n \in \N$ if necessary. This completes the
proof of (\ref{fmn}).

We have verified the first claim of the lemma, since it evidently follows from
(\ref{vyn}), subject to proving the lemma that follows. 

The second claim is proved analogously. The analog of (\ref{nmzerobd}) is
\begin{equation}\label{zerobd}
 \frac{\max_{1 \leq m \leq n} \Big( \sum_{l=1}^{n-m} q_l \Big)^2}{
   \sum_{m=1}^n  \Big( \sum_{l=1}^{n-m} q_l \Big)^2},
\end{equation}
which follows from the inequality
$$
 \Big( \sum_{l=1}^{n-1} q_l \Big)^2 \leq C n^{2 \alpha} L(n)^{-2},
$$
which is a consequence of Lemma \ref{lemprobr}(i), and the bound
$$
 \sum_{m=1}^n \Big( \sum_{l=1}^{n-m} q_l \Big)^2 \geq 
  c n^{2 \alpha + 1} L(n)^{-2},
$$
which follows from Lemma \ref{lemprobr}(i) and (\ref{uct}). $\Box$ 

It remains to show that the first hypothesis of Theorem \ref{martclt} holds.
Recalling (\ref{formneg}), we see that,
for $- k_n < M < 0$,
$\big\vert S_n^{M+1} - S_n^M \big\vert^2 \leq 4 \big( F_M^n \big)^2$,
whose right-hand-side is $o\big(\E V_n \big)$ by  
(\ref{fmn}), (\ref{vnrep}) and Lemma \ref{lemvn}. 
For
$0 \leq m < n$, we use (\ref{formpos}) to find that
$\big\vert S_n^{m+1} - S_n^m \big\vert^2 \leq 4 
 \big( \sum_{l=1}^{n-m} q_l \big)^2$. However, this right-hand-side is at
 most  $4 c n^{2\alpha} L(n)^{-2}$ by Lemma \ref{lemprobr}(i). From the
 discussion after (\ref{varsnf}), we know that $\E V_n \geq \big( 1-\epsilon_n
 \big) {\rm Var} S_n$ for $n$ sufficiently high, so that $4 c n^{2\alpha}
 L(n)^{-2} \leq O 
\big( n^{-1} \E(V_n) \big)$ by Lemma \ref{lemvarsn}.  
 Thus, the first hypothesis
 of Theorem \ref{martclt} is indeed satisfied. $\Box$
\end{section}
\begin{section}{The FKG inequality and convergence in $L^{\infty}$}
This section is devoted to the remaining step in the proof of Theorem \ref{thm}, namely, to improving the topology in our convergence result.
\begin{prop}\label{propfkg}
Let $\mu$ be any probability measure on $\mathbb{N}$, and let
$\lambda$ be any extremal $\mu$-Gibbs measure. Let $X,Y \in L^2 \big( \lambda \big)$ be two increasing functions, where
$\{ -1,1 \}^{\mathbb{Z}}$ is given the natural poset structure. Then
$X$ and $Y$ are not negatively correlated under $\lambda$.
\end{prop}

\noindent{\bf Proof.} Note that ${\rm Cov} \big( X,Y \big)
\geq 0$ if and only if
\begin{equation}\label{varxy}
{\rm Var} \big( X + Y \big) \geq {\rm Var} \big( X \big)
  + {\rm Var} \big( Y \big).
\end{equation}
For any $Z \in L^2 \big( \lambda \big)$, we write
\begin{equation}\label{eqz}
Z = \sum_{n \in \mathbb{Z}} \Big( \mathbb{E} \big( Z \big\vert \sigma_n \big) -  \mathbb{E} \big( Z \big\vert \sigma_{n-1} \big) \Big),
\end{equation}
where $\sigma_n$ is the $\sigma$-algebra generated by coordinates
with index strictly less than $n$. 
 In regard to (\ref{eqz}), note that the maps $L^2(\lambda) \to \R: Z \to \E(Z\vert \sigma_n)$ that project functions onto their mean given
$\sigma_n$ form an increasing sequence of projections in
the Hilbert space $L^2(\lambda)$. The identity (\ref{eqz}) expresses Z
as a sum of differences of the successive projections. As such, this
is a sum of orthogonal vectors in $L^2(\lambda)$. The sum of a countable collection of
orthogonal vectors in a Hilbert space has finite norm if and only if
the sum of the squares of the norms of these vectors is finite. The
norm of the sum has a square given by the sum of the squares of the
norms of the constitutent vectors. From this, we obtain
$$
{\rm Var}(Z) = \sum_{n \in \Z}{\E \big(  \mathbb{E} \big( Z \big\vert \sigma_n \big) -  \mathbb{E} \big( Z \big\vert \sigma_{n-1} \big) \big)^2}. 
$$
We find then that
\begin{equation}\label{varz}
 {\rm Var} (Z) = \sum_{n \in \mathbb{Z}} \mathbb{E}_{\sigma_{n-1}} {\rm Var}
  \Big(   \mathbb{E} \big( Z \big\vert \sigma_n \big) -  \mathbb{E} \big( Z \big\vert \sigma_{n-1} \big) \Big).
\end{equation}
Note that, in the summand on the right-hand-side of (\ref{varz}), ${\rm Var}$ denotes a conditional variance: the data in $\sigma_{n-1}$ is fixed, and the relevant randomness  arises from the bit with index $n-1$. 

%The proof of Proposition 4: page 22, line 5, (4.35): yes, in
%\E_{\sigma_{n-1}} Var (  E(Z \vert \sigma_n) - \E ( Z \vert
%\sigma_{n-1}) ),
%Var is a conditional variance: the data in \sigma_{n-1} is fixed, and
%the relevant randomness is bit n-1 (which we learn under \sigma_n but
%do not know under \sigma_{n-1}).
%To interpret the equation that precedes (4.35), note that the maps
%from L^2(\lambda) that project functions onto their mean given
%coordinates at most n form an increasing sequence of projections in
%the Hilbert space L^2(\lambda). The equation in questions expresses Z
%as a sum of differences of the successive projections. As such, this
%is a sum of orthogonal vectors. The sum of a countable collection of
%orthogonal vectors in a Hilbert space has finite norm if and only if
%the sum of the squares of the norms of these vectors is finite. The
%norm of the sum has a square given by the sum of the squares of the
%norms of the constitutent vectors. In our case, the squared norm of
%the generic term (indexed by n) is the variance of the term, which may
%be expressed as the mean of the conditional variance given data in
%slots less than n. Thus, (4.35) is a consequence of an infinite
%dimensional version of Pythagoras' theorem in Hilbert space.

Given $\sigma_{n-1}$, each of $X$ and $Y$ does not decrease if we condition on the value of the $n$-th coordinate to be $1$, and does not increase if this value is conditioned to be $-1$. Thus,
\begin{eqnarray}
 & & {\rm Var}
  \Big(   \mathbb{E} \big( X + Y \big\vert \sigma_n \big) -  \mathbb{E} \big( X + Y  \big\vert \sigma_{n-1} \big) \Big) \nonumber \\
 & \geq &  {\rm Var}
  \Big(   \mathbb{E} \big( X \big\vert \sigma_n \big) -  \mathbb{E} \big( X   \big\vert \sigma_{n-1} \big) \Big) + {\rm Var}
  \Big(   \mathbb{E} \big( Y \big\vert \sigma_n \big) -  \mathbb{E} \big( Y  \big\vert \sigma_{n-1} \big) \Big), \nonumber
\end{eqnarray}
$\sigma_{n-1}$-a.s. Taking expectation over $\sigma_{n-1}$ and summing yields (\ref{varxy}) by means of $(\ref{varz})$. $\Box$
\begin{definition}
Recall that $S_0 = 0$, $S_i =  \sum_{j=1}^i \walk_j$ if $i > 0$, and $S_i = - \sum_{j=i}^{-1} \walk_j$ if $i < 0$, where $\walk_i = \walk_n{\lambda_p}(i)$.
For $l < m$, write
$$
A_{l,m} = \max_{i \in \{ l,\ldots,m\} } \big( S_i - (2p-1)i \big) \, - \,
\big( S_l - (2p-1)l \big)
$$
and
$$
B_{l,m} = \max_{i \in \{ l,\ldots,m\} } \big( S_i - (2p-1)i \big) \, - \,
\big( S_m - (2p-1)m \big)
$$
\end{definition}
\begin{lemma}\label{lemmax}
There exists a constant $C_0 > 0$ such that, 
for each $\epsilon > 0$, there exists $n_0 = n_0(\epsilon)$
such that, for $n \geq n_0$,
$$
 \mathbb{P} \bigg( A_{0,n}  > \Big( C_0 + 50 \varcon^{-1} \sqrt{\log \big( \epsilon^{-1} \big)}
  \Big) 
 n^{1/2 + \alpha} L(n)^{-1} \bigg) < \epsilon,
$$
where the constant $\varcon = \varcon(p)$ was defined in (\ref{defvarcon}).
\end{lemma}
\noindent{\bf Proof.}
We write $A_n = A_{0,n}$
and $B_n = B_{0,n}$.
For  $n \in \mathbb{N}$,
set $h(n) = 3 \varcon^{-1} \sqrt{\log \big( \epsilon^{-1} \big)} n^{1/2 + \alpha}
L(n)^{-1}$. 
We begin by establishing the following statements. For any $\epsilon > 0$, there exists $n_0 = n_0(\epsilon)$ such that, if $n \geq n_0(\epsilon)$, then, for any $K > 0$,
\begin{equation}\label{ban}
\mathbb{P}\Big(  B_n  > K + h(n) \Big) \geq \epsilon  \quad \implies \quad
\mathbb{P}\big(  A_n > K \big) \geq 1 - \epsilon
\end{equation}
and
\begin{equation}\label{abn}
\mathbb{P}\big(  A_n > K + h(n) \big) \geq \epsilon  \quad \implies \quad
\mathbb{P}\Big(  B_n  > K \Big) \geq 1 - \epsilon.
\end{equation}
Note that $A_n$ is an increasing and $B_n$ a decreasing random variable.
As such, supposing that $\mathbb{P} \big( B_n  > K + h(n)  \big) \geq \epsilon$ and  $\mathbb{P} \big( A_n \leq K   \big) \geq \epsilon$,
then
\begin{eqnarray}
 & & \mathbb{P} \Big( \Big\{ B_n  > K + h(n) \Big\} \cap \Big\{  A_n \leq K \Big\} \Big) \nonumber \\
 & = & \mathbb{P} \Big( B_n  > K + h(n) \Big) \mathbb{P} \Big(
 A_n \leq K \Big\vert B_n  > K + h(n) \Big) \nonumber \\
 & \geq &  \mathbb{P} \Big( B_n  > K + h(n) \Big)  \mathbb{P} \Big(   A_n \leq K \Big) \geq \epsilon^2, \nonumber
\end{eqnarray}
the first inequality by means of Proposition \ref{propfkg}.
However,
$A_n - B_n  = S_n - (2p-1)n$, so that
$$
\mathbb{P} \Big(  B_n - A_n > h(n) \Big)
 = \mathbb{P} \Big( n^{-1/2 - \alpha} L(n) \varcon \big( S_n - (2p-1)n \big) 
           < - 3 \sqrt{\log \big(
   \epsilon^{-1} \big)} \Big) < \epsilon^2,
$$
the inequality valid for $n$ sufficiently high, because Proposition
\ref{propfd}
implies that $n^{-1/2 - \alpha} L(n) \varcon \big( S_n - (2p-1)n \big)$ converges in distribution to a
Gaussian random variable of mean zero and unit variance. 
In this way, we establish 
(\ref{ban}), with (\ref{abn}) following similarly.

We now show that there exists $\alpha_n \in (0,\infty)$ such that, for $\epsilon > 0$ and $n \geq n_0(\epsilon)$,
\begin{equation}\label{anint}
\mathbb{P} \bigg(  A_n \in \Big(  \alpha_n - 3h(n) , \alpha_n + 3h(n)  \Big) \bigg) \geq  1 - 2\epsilon
\end{equation}
and
\begin{equation}\label{bnint}
\mathbb{P} \bigg(  B_n \in \Big(  \alpha_n - 3h(n) , \alpha_n + 3h(n)  \Big) \bigg) \geq  1 - 2\epsilon.
\end{equation}
To this end, set 
\begin{equation}\label{eqkn}
K(n) = \sup \big\{ k \in (0,\infty): \mathbb{P} \big( A_n > k \big) \geq \epsilon  \big\} -\epsilon/2.
\end{equation}
Note that $\mathbb{P} \big( A_n > K(n) \big) \geq \epsilon$ and
$\mathbb{P} \big( A_n > K(n) + \epsilon \big) < \epsilon$. The first inequality forces $\mathbb{P} \big( B_n > K(n) - h(n) \big) \geq 1 - \epsilon$ by means of (\ref{abn}), which gives 
\begin{equation}\label{eqintcalc}
\mathbb{P} \big( A_n > K(n) - 2 h(n) \big) \geq 1 - \epsilon. 
\end{equation}
by  (\ref{ban}). We now set $\alpha_n = K(n) - h(n)$
and note that, if $\epsilon < h(n)$, then (\ref{anint}) holds.
Indeed, these choices for $\alpha_n$ and $\epsilon$, alongside  $P(A_n > K(n) + \epsilon ) < \epsilon$, yield
\begin{equation}\label{eqalpnan}
P\big( A_n > \alpha_n + 2h(n) \big) < \epsilon.
\end{equation}
The bounds (\ref{eqintcalc}) and (\ref{eqalpnan}) imply that
$$
P  (   \alpha_n + 2h(n) \geq A_n \geq \alpha_n - h(n)   )   \geq 1 - 2\eps,
$$
whence (\ref{anint}).
Note now that
$\mathbb{P} \big( B_n \geq \alpha_n + 2h(n) + \epsilon \big) < \epsilon$: for otherwise, (\ref{ban}) would imply that
$\mathbb{P} \big( A_n > \alpha_n + h(n) + \epsilon \big) \geq 1 - \epsilon$, contradicting our assumption. We see that
$$
\mathbb{P} \Big( B_n \in \big( \alpha_n , \alpha_n + 2h(n) + \epsilon \big) \Big) \geq 1 - 2 \epsilon,
$$
whence (\ref{bnint}).

We set $r_n = \alpha_n n^{-1/2 - \alpha} L(n)$. We claim that, 
for any $K \in \mathbb{N}$, there exists $n_0 \in \mathbb{N}$ 
such that $n \geq n_0$ implies
\begin{equation}\label{rkn}
 r_{Kn} \leq 2 K^{-1/2 - \alpha} r_n
 +  20 \varcon^{-1} \sqrt{\log \big( \epsilon^{-1} \big)} 
 +  24 \varcon^{-1} \sqrt{\log \big( \epsilon^{-1} \big)} K^{-1/2 - \alpha}.
\end{equation}

%"we readily conclude" on line 1 of page 24: in fact, we do not take K
%to infinity, but fix it at some large value (such that $2 K^{-1/2 -
%\alpha} < 1$).
%(4.40) then forces the limsup of r_m over the subsequence $m = Kn$ to
%be at most $(20 + 24 + 1) \tilde{c}^{-1}\sqrt{}\log(\epsilon&{-1})$.
%The restriction that the limsup be taken over this subsequence can
%then be dropped, since \vert \alpha_m - \alpha_n \vert = \Theta(\vert
%m - n \vert), since $A_n$ and $B_n$ are processes whose increments are
%all at most one.

From this, we will argue that
\begin{equation}\label{eqreadconc}
\sup_{n \in \mathbb{N}} r_n \leq C_0 +  40 \varcon^{-1} \sqrt{\log \big( \epsilon^{-1} \big)} 
\end{equation}
for some $C_0 > 0$. Indeed, fixing $\kmac \in \N$ in (\ref{rkn}) high enough that
$2 \kmac^{-1/2 - \alpha} < 1$, we find that the limsup of $r_m$ along the subsequence 
$\big\{ m = \kmac n: n \in \N \big\}$ is at most $(20 + 12 + 1) \big( \varcon^{-1} \sqrt{\log \big( \epsilon^{-1} \big)} \big)$. Note then that $\big\{ A_n: n \in \N \big\}$ is a random process whose increments are in absolute value at most one. Thus, the definition (\ref{eqkn}) implies that $\vert \alpha_m - \alpha_n \vert \leq \vert m - n \vert$. From this, we obtain the above bound on the limsup of $r_m$ over all $m \in \N$ (with an arbitrarily small addition to the value $33$); and thus (\ref{eqreadconc}).

In tandem with (\ref{anint}), (\ref{eqreadconc}) yields
$$
\mathbb{P} \bigg( A_n > \Big( C_0 + 49 \varcon^{-1}  \sqrt{\log \big( \epsilon^{-1} \big)} \Big)  n^{1/2 + \alpha} L(n)^{-1} \bigg) \leq 2\epsilon,
$$
from which the statement of the lemma follows.
To derive (\ref{rkn}), note that
\begin{equation}\label{akn}
A_{Kn} \leq \max_{j \in \{ 1,\ldots,K-1\} } \Big( S_{jn} - \big( 2p - 1 \big) j n  \Big) 
 \, + \, \max_{j \in \{ 0,\ldots, K - 1\} } A_{jn,(j+1)n}.
\end{equation}
By the convergence for finite-dimensional distributions stated in
Proposition \ref{propfd}, 
we have that, for each $\epsilon > 0$, there exists
$n_0 = n_0(\epsilon)$ such that, for $n \geq n_0$,
$$
\mathbb{P} \Big(  \max_{ j \in \{ 1,\ldots,K-1\} } 
 \Big( S_{jn} - \big( 2p - 1 \big) j n  \Big)   >   
 3 \varcon^{-1}  \sqrt{\log \big( \epsilon^{-1} \big)}  n^{1/2 +
  \alpha}  L(n)^{-1} \Big) < \big( K - 1 \big) \epsilon.
$$
Note also that
$$
\mathbb{P} \Big( \max_{j \in \{ 0,\ldots, K - 1\} } A_{jn,(j+1)n} \geq \alpha_n + 3h(n) \Big) \leq 2K \epsilon,
$$
by (\ref{anint}). From (\ref{akn}), then,
\begin{equation}\label{thrkeps}
\mathbb{P} \Big( A_{Kn} > 
 3 \varcon^{-1}  \sqrt{\log \big( \epsilon^{-1} \big)}  n^{1/2 + \alpha} L(n)^{-1} + \alpha_n + 3 h(n)\Big) \leq \big( 3K - 1 \big) \epsilon
\end{equation}
for $n \geq n_0(\epsilon)$.
Choosing $\epsilon > 0$ to satisfy 
$(3K-1) \epsilon \leq 1 - 2 \epsilon$, from (\ref{thrkeps})
and the inequality arising from (\ref{anint}) by the substitution of $Kn$ for $n$, we see that
$$
 3 \varcon^{-1}  \sqrt{\log \big( \epsilon^{-1} \big)}  n^{1/2 + \alpha} L(n)^{-1} 
 + \alpha_n + 3 h(n) \geq \alpha_{Kn} - 3h \big( Kn \big)
$$
Rearranging, and by $L:(0,\infty) \to (0,\infty)$ being a slowly varying
function, we obtain (\ref{rkn}). $\Box$ \\
\noindent{\bf Proof of Theorem \ref{thm}.} 
We will make use of:
\begin{lemma}\label{lemsupapp}
There exists a constant $C \in (0,\infty)$ such that, for any
$\epsilon_0 > 0$ sufficiently small and for $\delta > 0$,
$$
\limsup_{n \to \infty} \mathbb{P} \Big( \sup_{k = 0,\ldots,\lfloor \epsilon_0^{-1}
  \rfloor} A_{k \lfloor \epsilon_0 n \rfloor, (k+1) \lfloor \epsilon_0 n
  \rfloor} > \delta n^{1/2 + \alpha}L(n)^{-1} \Big) 
\leq  C \big( \epsilon_0^{-1} + 1 \big) \exp \Big\{ - \frac{\varcon^2}{2\cdot 50^2} \delta^2 
 \epsilon_0^{-1/2 - 2\alpha} \Big\}.
$$
\end{lemma}
\noindent{\bf Remark.} An analogous bound holds for the counterpart of $A_{l,m}$, in which the maximum taken in the definition is replaced by a minimum.\\
\noindent{\bf Proof.} Note that (\ref{lureps}) implies that, for
$n$ sufficiently high,
$$
 n^{1/2 + \alpha} L(n)^{-1} 
  \geq \epsilon_0^{-1/4 - \alpha} 
  \big( \epsilon_0 n\big)^{1/2 + \alpha} L \big( \epsilon_0 n \big)^{-1}.
$$
The probability that we must estimate is bounded above by
$$
\big( \epsilon_0^{-1} + 1 \big) \mathbb{P} \Big( A_{0,\lfloor \epsilon_0 n
  \rfloor} \geq \delta \epsilon_0^{-1/4 - \alpha} \big( \epsilon_0 n
\big)^{1/2 + \alpha} L \big( n \epsilon_0 \big)^{-1} \Big) 
$$
Setting $\epsilon > 0$ in Lemma \ref{lemmax} according to
$
C_0 + 50 \varcon^{-1} \sqrt{\log \big( \epsilon^{-1} \big)} = \delta 
 \epsilon_0^{-1/4 - \alpha}$, we obtain that the last displayed expression
 is at most
$$
 \big( \epsilon_0^{-1} + 1 \big) \exp \Big\{ -  \frac{\varcon^2}{50^2} \Big( \delta 
 \epsilon_0^{-1/4 - \alpha} - C_0 \Big)^2 \Big\}
\leq  C  \big( \epsilon_0^{-1} + 1 \big) \exp \Big\{ -  \frac{\varcon^2}{2 \cdot 50^2} \delta^2 
 \epsilon_0^{-1/2 - 2 \alpha}  \Big\}
$$ 
for some constant $C>0$, for $\epsilon_0$ small, and for $n$ sufficiently high. $\Box$ \\
For the proof of the theorem, it suffices to construct, for each $\epsilon
> 0$, a sequence $\big\{ {\rm C}_n^\epsilon: n \in \N  \big\}$ of couplings
$S_p^n$ and $\fbm_{\alpha + 1/2}$ for which
\begin{equation}\label{coupcne}
 \lim_{n \to \infty} {\rm C}_n^\epsilon 
 \Big(   \vert\vert  S_p^n - \fbm_{\alpha + 1/2} \vert\vert_{L_{\infty}\big( [0,T] \big)} > \epsilon  \Big)  = 0,
\end{equation}
for then we may find an increasing sequence $\big\{ n_i: i \in \N \big\}$ 
such that
$$
 \sup_{n \geq n_i} {\rm C}_n^{2^{-i}} 
 \Big(   \vert\vert  S_p^n - \fbm_{\alpha + 1/2} \vert\vert_{L_{\infty}\big(
   [0,T] \big)} > 2^{-i}  \Big) < 2^{-i},
$$
for each $i \in \N$, and set ${\rm C}_m = {\rm C}_m^{2^{-k_m}}$,
where $k_m$ is the maximal $n_i$ that does not exceed $m$; note that $k_m
\to \infty$ ensures that this sequence of couplings is as the statement of
the theorem demands.
A coupling ${\rm C}_n^{\epsilon}$ suitable for (\ref{coupcne}) may be obtained by the
use of Proposition \ref{propfd} 
to couple the values of $S_p^n$ and $\fbm_{\alpha + 1/2}$ at points of the
form 
$\big\{ i\epsilon T: i = 0,\ldots, \lfloor \epsilon^{-1} \rfloor + 1
\big\}$, so that the maximum difference $\vert S_p^n - \fbm_{\alpha +1/2} \vert$ at
such points tends to zero in probability under ${\rm C}_n^{\epsilon}$. 
We
then use Lemma \ref{lemsupapp}, and the remark following its statement,  as well as $\fbm_{\alpha + 1/2}$ being uniformly continuous on
$[0,T]$ to verify (\ref{coupcne}). $\Box$
\end{section}

\begin{section}{Open problems}
\begin{itemize}
\item
It would be interesting to find the negatively
correlated fractional Brownian motions, with Hurst parameter
$H\in(0,{1/2})$, as scaling limits of variants of the discrete
processes that we consider. A natural first guess is the walk arising from
the model in which each vertex of $\Z$ is given the opposite sign of
its parent in $G_\mu$, instead of the same sign.
 However, it is quite possible that this gives rise
to Brownian motion as a scaling limit of the associated walk. 
A more promising candidate discrete model is one in which the vertices
in each component of $G_\mu$ --- ordered with respect to the natural
ordering of $\mathbb Z$ --- are given values that alternate between $1$ and
$-1$, independently choosing one of the two ways
of doing this, each with probability $1/2$.
\item As mentioned in the discussion following the statement of Theorem
  \ref{thm}, we believe that the theorem is sharp. We pose the problem to 
show that if a measure $\mu$ on $\N$ is such that the conclusion of
  Theorem \ref{thm} holds, with some deterministic function of $n$ playing
  the role of $\tilde{c} n^{-1/2 - \alpha}L(n)$, then $\mu \in
  \Gamma_\alpha$. An important step here would be to show that, if the measure
  $\mu$ is such that the variance of $S_n$ is a slowly varying multiple of
  $n^{2\alpha + 1}$, then $\mu \in \Gamma_\alpha$. 
\item  Define $\tilde{\Gamma}_{\alpha}$  
according to Definition \ref{defone}, with the alteration that $\mu$ is supported
on $\Z$ and is symmetric about $0$, and let $\mu
\in \tilde{\Gamma}_{\alpha}$ for some $\alpha \in (0,\infty)$. Consider a voter
model in dimension one, in which voters reside at the elements of
$\Z$. Each voter at a given time has an affiliation to one of two political
parties. At any given moment $t \in \Z$ of time, each
voter selects the resident at a displacement given by an independent sample
of $\mu$, and inherits the affiliation that this resident held at time
$t-1$. (A continuous time variant, where opinions are imposed on a voter at
a Poisson point process of times, may also be considered.) 
We propose the problem of considering the set of equilibrium measures of
this process as an analogue of the family $\big\{\lambda_p : p \in [0,1]
\big\}$ and seeking a counterpart to Theorem \ref{thm}. The two unanimous
configurations are always equilibrium measures, and we anticipate that
mixtures of these are the only such measures exactly when $\alpha \geq
1$ (and the condition ${\rm g.c.d.}(\mu) = 1$ is satisfied). 
This is because the transmission histories of the affiliation held by
two voters will almost surely coincide in the distant past precisely when a 
discrete time random walk with step distribution $\mu$ almost surely
reaches zero. The case $\alpha = 1$ corresponds to the step
distribution of the $x$-displacement of a two-dimensional simple random
walk between its successive visits to the $x$-axis; such a walk is
recurrent, but only ``marginally so'', pointing to the value of $\alpha =
1$ being critical for the problem.

\item Let $d \geq 2$. 
For any measure $\mu$ whose support is contained in the integer lattice 
$\mathbb{Z}^d$, the random $\Z^d$-spanning graph structure $G_{\mu}$ may be
defined.  Consider a law $\mu$ that has a regularly decaying heavy tail. 
For example, we might
  insist that 
$\mu \big( B_n^c \big) =  n^{-\alpha} L\big(\vert\vert n \vert\vert \big)$,
where $L:(0,\infty) \to (0,\infty)$ is slowly varying. Here, $B_n = \big\{
x \in \Z^d: \vert\vert x \vert\vert \leq n \big\}$.   
If $\mu$ is chosen to be ``symmetric'', in the sense that $\mu \{ x \} =
\mu \{ y \}$
whenever $x = y$, it is not hard to show that all the components of
$G_{\mu}$ are finite, each containing a unique cycle. Choices of $\mu$
supported in a half-plane $\big\{ x_1 < 0\big\}$ will give rise to
infinite components, however. This raises the question of finding a phase
transition in $\alpha$ for uniqueness of the infinite component for some
family of measures $\mu$ supported in the strict half-plane  
 $\big\{ x_1 < 0\big\}$. In fact, the discussion in the preceding problem
 already addresses this problem in a certain guise. To obtain the voter model problem,
 we take $d = 2$, write points
 in the plane in the form $(t,x)$, take $\mu$ to be supported on the line
 $\big\{ x  = -1 \big\}$, with $\mu\{ -1, \cdot \} \in \tilde{\Gamma}_{\alpha}$.
Naturally, we might take $d > 2$. In this case, in seeking a result
analogous to Theorem \ref{thm}, we would seek to show that, for regularly
varying laws $\mu$ whose decay is slow enough to ensure that the transition
histories of distinct voters may be disjoint, the affiliation ``white
noise'' in the time-slice $\big\{ t = {\rm constant} \big\}$ at equilibrium is given by a
fractional Gaussian field.
\item 
Theorem \ref{thm} states that fractional simple random walk and
fractional Brownian motion may be coupled to be close in $L^{\infty}$ on
compact sets under rescaling of the discrete walk. One might seek to
quantify this, investigating how quickly $\epsilon$ may be taken to $0$
in the limit of high $n$ in (\ref{cncoupl}), for some coupling $C_n$. The analogue
for simple random walk and its convergence to Brownian motion is the
celebrated Komlos-Major-Tusnady theorem \cite{kmt}, one case of which
states that 
simple random walk may be coupled with asymptotic probability one to a Brownian motion with a uniform
error in the first $n$ steps of at most a large constant multiple of $\log n$. 
Quantifying the rate of convergence
could be useful if we would like to compare the stochastic difference equation
driven by the fractional discrete process and its continuous analogue
driven by fractional Brownian motion, the latter having been extensively
studied. (See \cite{bookone,dhp,booktwo} for treatments of 
the stochastic calculus of fractional Brownian motion.)
\end{itemize}
\end{section}

\appendix

\begin{section}{Market interpretation}
This section was added in response to referee's question about whether there was a market-based ``story'' to go with our random walk model.  We will not address the empirical question of whether our construction describes actual markets.  The literature on market-increment autocorrelations is too vast (searches for {\em market momentum} or {\em mean reversion} turn up hundreds or thousands of articles per year) to summarize here.  Instead, we offer a few very general observations in response to a simpler question: {\em if} an asset's price were described by one of our random walks (for some $\mu$), what kind of story would explain it? Could one reconcile the story with (some form of) market efficiency?   The first story to spring to mind is the following:

\noindent{\bf Story 1: Peer-influenced decisions.} {\it Many decisions involve the consultation (conscious or otherwise) of people who
made a similar decision recently.  Jeans or khakis?  Buy or rent?  State or private school?
If the decision probability is a monotone linear function of previous decisions, one can represent this by having individuals copy (with some probability) a randomly chosen previous decision.  The time elapsed since that decision is random (and a power law is not unreasonable).  The cumulative number of decisions of one type, minus the number of the other type, is thus one of the random walks described in this paper.}

\vspace{0.2cm}

This story might explain why some of a firm's fundamentals (market share, revenue, etc.) would correspond to our model for some $\mu$.  It does not explain momentum effects in asset prices, however, since one should be able to foresee these effects and price them in.

In classical finance, asset prices are martingales with respect to the {\em risk neutral} probability (assuming no interest, a point we return to later).  This may in fact be taken as the definition of risk neutral probability.  In models with $\pm 1$ increments, the risk neutral probability measure is that of the simple random walk.  (Let us assume that the probability of an up step, conditioned on the past, is always strictly between zero and one.)  Assuming no arbitrage, the price of a derivative is its expected value in this measure.
Because of this simplicity, models with $\pm 1$ increments are especially natural to work with.\footnote{Any continuous model $X_t$ can be interpreted as having $\pm 1$ increments if one changes time parameterization: simply define times $t_k$ such that $X_{t_k}$ is an integer for each $k$ and $t_{k+1}$ is the first time after $t_k$ that a distinct integer is reached.  Then $N_k := X_{t_k}$ has $\pm 1$ increments.  One can even approximate a discontinuous jump process by using a probability measure in which $N_k$ is likely to sometimes go up or down several steps in a row.}

The interesting question is the following: {\em why do risk neutral and true probability differ?}

\vspace{.1in} 

{\em \bf Note:} Discrepancies between true probability and risk neutral probability do not necessarily imply market inefficiency.  However, they are most plausible if the discrepancy is small.  (If a trader could get a $5000$ percent return with probability $.99$ in one year, this would be hard to reconcile with market efficiency.)  In our model, this suggests using a $\mu$ that decays slowly.
In the classical capital asset pricing model (CAPM) these discrepancies are explained by the asset's correlation with so-called {\em systemic risk}: when the price goes down, the so-called {\em market portfolio} is generally likely to go down also; the demand for money in that scenario is greater because people are risk averse.  In the simplest CAPM story, all investors hold some combination of risk free assets and the market portfolio.

\vspace{0.2cm}

\noindent{\bf Story 2:  Inhomogeneous market portfolio.} {\it Suppose that one population of traders tends to enter and leave the market for apparently irrational reasons.  Another population of ``savvy traders'' holds a time inhomogenous ``savvy market portfolio,'' including long and short positions that vary in time.  This group makes money in expectation by holding this portfolio (twelve percent per year, say), while the other traders make less money in expectation.}

\vspace{0.2cm}

The behavior of the savvy traders in the above story makes perfect sense, but we have not explained the actions of the ``less savvy'' group.  Let us suppose that the asset price is $10$ today and has gone up recently, so that it now has a $.52$ chance to be $11$ tomorrow, and a $.48$ chance to be $9$ tomorrow.  Or suppose the asset price has dropped recently and now has a $.48$ to chance to be $11$ tomorrow and a $.52$ to be $9$.  In each case, it is clear why a savvy trader would take the side of the bet with positive expectation.  Under certain parameters, it is also reasonable to suppose that such investors are too risk averse or otherwise constrained to bid the price all the way up or down to tomorrow's expected value.  But who is taking the negative-expectation side of the bet?  One can think of many stories, but here is one that fits our model.

\vspace{0.2cm}

\noindent{\bf Story 3:  Anti-momentum trading.} {\it Less savvy traders tend to sell after the price goes up (``profit taking") and buy after the price goes down (``bargain hunting").  In other words, they systematically bet against the momentum strategy.  Each less savvy trader decides whether to hold a stock during a period of time by looking at price changes in the recent past.  The probability that the less savvy trader takes the ``less savvy'' position is a linear function of these prior changes.  The amount this skews the price is approximately linear in the number of less savvy traders taking the anti-momentum strategy.}

\vspace{0.2cm}

This story suggests a rather paradoxical conclusion.  Market returns may exhibit positive autocorrelation (momentum) precisely because the less savvy investors {\em think} or {\em feel} that they should exhibit negative autocorrelation (mean reversion).  Indeed, the fact that the savvy investors tend to earn higher returns essentially {\em implies} that whatever strategy the less savvy investors employ tends to be wrong.

We note a couple of other (possibly complementary, possibly second-order) stories:

\vspace{0.2cm}

\noindent{\bf Story 4: Variable interest rate.}  {\it The presence of interest affects the definition of risk neutral probability.  The interest-discounted asset price (not the asset price itself) is the martingale.  Even if the rate of interest earned per time unit is constant or very slowly varying, it may become non-constant when time is parameterized by the number $k$ of integer price changes.  Perhaps when the price has gone up recently, people trade more slowly (so
that the interest per tick is higher). This leads to momentum affects when one parameterizes time by $k$.}

\vspace{0.2cm}

\noindent{\bf Story 5:  Non-liquidity/inefficiency.}
{\it There {\em is} inefficiency, but the market is too small for the arbitrage opportunities to be very valuable.  There may
be times at which price history indicates a high likelihood of a rise or fall in future prices, but (for whatever reason)
the volume of trading at these times is not sufficient to attract arbitrageurs.}

\vspace{0.2cm}

Finally, we remark that if a market were to exhibit momentum on the scale of days and mean reversion on the scale of several months, this could be modeled with a variant of our walks, where one samples $j$ according to $\mu$ and copies the increment $j$ steps previously (for some range of $j$ values) or the opposite of that increment (for another range of $j$ values).  The stories described above make sense for these variants as well.
However, we stress, in conclusion, that the extent to which any of the above stories is, in actual markets, {\em correct} is an empirical matter well beyond the scope of this work.
\end{section}

\bibliographystyle{plain}
%\bibliography{fbmbib}

\end{document}